\documentclass[12pt,a4paper]{article}
\usepackage{amsfonts}
\usepackage{amsfonts,amssymb,amsmath,indentfirst,amsthm}

\input{tcilatex}

\begin{document}

\date{}
\title{Regularity for weak solutions to nondiagonal quasilinear degenerate elliptic systems}
\author{Yan Dong, Pengcheng Niu \thanks{%
Pengcheng Niu: Corresponding author, Department of Applied
Mathematics, Northwestern Polytechnical University, Xi'an, Shaanxi,
710129, China. e-mail: pengchengniu@nwpu.edu.cn}
}
\date{}
\maketitle

\newtheorem{thm}{{\indent}Theorem} \newtheorem{cor}[thm]{{\indent}Corollary}
\newtheorem{lem}[thm]{{\indent}Lemma} \newtheorem{prop}[thm]{%
{\indent}Proposition} \theoremstyle{definition} \newtheorem{defn}[thm]{%
{\indent}Definition}

\newtheorem{theo}{\hspace*{2.0em}Theorem\hspace*{0.1em}}[section]
\newenvironment{keywords}{\par\textbf{keywords:}\mbox{ }}{ } \newenvironment{%
coj}{\par\textbf{Conjecture:}\mbox{ }}{ } \numberwithin{equation}{section}

\textbf{Abstract.} The aim of this paper is to establish regularity
for weak solutions to the nondiagonal quasilinear degenerate
elliptic systems related to H\"{o}rmander's vector fields, where the
coefficients are bounded with vanishing mean oscillation. We first
prove $L^p$($p \ge 2$) estimates for gradients of weak solutions by
using a priori estimates and a known reverse H\"{o}lder inequality,
and consider regularity to the corresponding nondiagonal homogeneous
degenerate elliptic systems. Then we get higher Morrey and Campanato
estimates for gradients of weak solutions to original systems and
H\"{o}lder estimates for weak solutions.




\section{Introduction}

\label{introduction}

Regularity for solutions to elliptic systems in Euclidean spaces has
been studied by many authors and a lot of important conclusions were
got. Campanato in [2] obtained gradient estimates for weak solutions
to linear elliptic system with discontinuous coefficients. For
related articles, we quote [1, 14] and the references therein.

Huang in [18] derived Morrey estimates for uniformly elliptic
systems by applying Campanato's technique. Zheng and Feng in [28]
established H\"{o}lder estimates for weak solutions to quasilinear
elliptic systems by reverse H\"{o}lder inequality and Dirichlet
growth theorem, where the coefficients belong to $L^\infty  \cap
VMO$, and low terms satisfy controlled growth conditions. For the
following second order quasilinear elliptic systems
\[
 - D_\alpha  a_i^\alpha  (x,u,Du) = a_i (x,u,Du),
\]
where $a_i^\alpha  (x,u,Du) = A_{ij}^{\alpha \beta } (x)D_\beta u^j
+ g_i^\alpha  (x,u,Du)$, $A_{ij}^{\alpha \beta }  \in C^{0,\alpha
}$, Dan\v{e}\v{c}ek in [4] proved Morrey and Campanato estimates
with $p=2$ for weak solutions. When the coefficients $A_{ij}^{\alpha
\beta }$ belongs to $L^\infty  (\Omega ) \cap \mathcal {L}_\phi
(\Omega )$ (where $\phi  = \frac{1} {{1 + \left| {\ln r} \right|}}
$), $A_{ij}^{\alpha \beta }$ belongs to $VMO(\Omega ) \cap L^\infty
(\Omega )$, or $A_{ij}^{\alpha \beta }$ is bounded and belongs to
Campanato spaces, Dan\v{e}\v{c}ek and Viszus in [5, 6, 7] gave
similar estimates. Chiarenza, Franciosi and Frasca ([3]) obtained
$L^p$ estimate for weak solutions to divergence linear elliptic
systems by representation formula.

To nondiagonal elliptic systems, Kawohl in [19] investigated
H\"{o}lder continuity for bounded weak solutions to qualilinear
elliptic systems if the Liouville type property for these systems is
satisfied. Wiegner in [25] gained H\"{o}lder regularity for weak
solutions to nondiagonal systems with low terms satisfying natural
conditions. More related results also see [10, 15, 21, 22, 24, 30]
and the references therein.

Regularity of degenerate elliptic systems formed by H\"{o}rmander's
vector fields ([17]) has received wide attention in recent years. Di
Fazio and Fanciullo in [8] proved Morrey estimates ($p=2$) for weak
solutions to linear degenerate elliptic systems. Dong and Niu in [9]
showed Morrey estimates ($ p \ge 2 $) for linear degenerate elliptic
systems. For nonlinear systems, Xu and Zuily in [26] handled
interior regularity of weak solutions to quasilinear degenerate
elliptic systems with the low term satisfying the natural condition.
Gao, Niu and Wang in [11] settled partial H\"{o}lder regularity for
weak solutions to degenerate quasilinear elliptic systems with the
coefficients belonging to $ VMO \cap L^\infty$ and the low term
satisfying the natural condition. We mention that those systems in
[8, 9, 11, 26] are all diagonal.

To our knowledge, there is no any regularity result to nondiagonal
degenerate elliptic systems. Whether do they have regularity? What
is the kind of regularity if they have? These are what we will
answer in this paper. Concretely, we consider the following
nondiagonal quasilinear degenerate elliptic system
\begin{equation}
\label{eq1} - X_\alpha ^ *  (a_{ij}^{\alpha \beta } (x,u)X_\beta u^j
) = g_i (x,u,Xu) - X_\alpha ^ *  f_i^\alpha  (x,u,Xu) ,
\end{equation}
where $X_\alpha   = \sum\limits_{k = 1}^n {b_{\alpha k}
(x)\frac{\partial } {{\partial x_k }}}$ ($b_{\alpha k} (x) \in
C^\infty  (\Omega )$) are real smooth vector fields in a
neighborthood $\tilde \Omega$ of some bounded domain $\Omega \subset
\mathbb{R}^n$ ($q \le n$) and satisfy H\"{o}rmander's condition of
step $s$ (see Section 2), $\alpha ,\beta = 1,2, \ldots ,q ;i,j =
1,2, \ldots , N; X_\alpha ^ *   =  - X_\alpha   + c_\alpha (c_\alpha
=  - \sum\limits_{k = 1}^n {\frac{{\partial b_{\alpha k} }}
{{\partial x_k }}}  \in C^\infty  \left( \Omega \right))$ is the
transposed vector field of $ X_\alpha$.

The aim is to establish higher integrability of gradients for weak
solutions to (\ref{eq1}), higher Morrey estimates, H\"{o}lder
estimates and higher Campanato estimate.

Before stating our main results, we need several assumptions of
(\ref{eq1}) (the detailed description for notions of Sobolev space
$W_X^{k,p}$, Morrey space $ L_X^{p,\lambda }$, Campanato space
$\mathcal {L}_X^{p,\lambda } $, H\"{o}lder space $C_X^{0,\kappa } $,
$BMO_X $ and $VMO_X $ sees Section 2):

(H1) Let $a_{ij}^{\alpha \beta } (x,u) = A^{\alpha \beta } (x)\delta
_{ij}  + B_{ij}^{\alpha \beta } (x,u)$, where $A^{\alpha \beta } (x)
\in VMO_X  \cap L^\infty, A^{\alpha \beta } (x)=A^{\beta \alpha }
(x), A^{\alpha \beta } (x)$ satisfy the ellipticity condition and
$B_{ij}^{\alpha \beta } (x,u)$ are bounded and measurable, that is,
there exist positive constants $ \lambda _0 , \mu _0 , \delta,  0 <
\lambda _0  \le \mu _0 , 0 < \delta  < 1,$ such that for $ a.e.$ $x
\in \Omega $ and for any $\xi  \in \mathbb{R}^{qN}$,
\[
\mathop {\lim }\limits_{R \to 0} \eta _R \left( {A^{\alpha \beta }
(x)} \right) = 0,
\]
\[
\lambda _0 \left| \xi  \right|^2  \le A^{\alpha \beta } (x)\xi
_\alpha  \xi _\beta   \le \mu _0 \left| \xi  \right|^2,
\]
\[
\left| {B_{ij}^{\alpha \beta } (x,u)} \right| \le \delta \lambda _0.
\]

(H2) Let $u \in W_X^{1,2} (\Omega ,\mathbb{R}^N ), g_i (x,u,z)$ and
$f_i^\alpha  (x,u,z)$ satisfy
\[
\left| {g_i (x,u,z)} \right| \le g^i (x) + L\left| z \right|^{\gamma
_0 } ,
\]
\[
\left| {f_i^\alpha  (x,u,z)} \right| \le g_i^\alpha  (x) + L\left| z
\right|,
\]
\[
f_i^\alpha  (x,u,z)z_\alpha ^i  \ge \gamma _1 \left| z \right|^2  -
\left( {g(x)} \right)^2 ,
\]
where $L$ and $\gamma _1$ are positive constants, $z \in
\mathbb{R}^{qN}$, $g^i  \in L_X^{pq_0 ,\lambda q_0 } (\Omega ), 1
\le \gamma _0 < \frac{{Q + 2}} {Q}$, $g_i^\alpha $ and $g \in
L_X^{p,\lambda } (\Omega )$, $p \ge 2, 0 < \lambda < Q, q_0 =
\frac{Q} {{Q + 2}}$,  $Q$ is the local homogeneous dimension
relative to $\Omega$ (see Section 2). Afterwards, we briefly denote
$ \tilde g = (g^i ), \tilde {\tilde g }= (g_i^\alpha).$

If $u \in W_X^{1,2} (\Omega ,\mathbb{R}^N )$ and for any $\varphi
\in C_0^\infty  (\Omega ,\mathbb{R}^N )$,
\[
 - \int_\Omega  {a_{ij}^{\alpha \beta } (x,u)X_\beta  u^j X_\alpha  \varphi ^j dx}
 = \int_\Omega  {\left( {g_i (x,u,Xu)\varphi ^i  - f_i^\alpha  (x,u,Xu)X_\alpha  \varphi ^i }
 \right)dx},
\]
we say that $u$ is a weak solution to (\ref{eq1}).

The main results of this paper are the following.

\textbf{Theorem 1.1} (higher integrability of gradients for weak
solutions) Let $u \in W_X^{1,2} (\Omega ,\mathbb{R}^N )$ be a weak
solution to (\ref{eq1}), the coefficients $a_{ij}^{\alpha \beta }$
satisfy (H1), $g_i$ and $f_i^\alpha$ satisfy (H2). Then there exists
a positive constant $\varepsilon _0  > 0$ such that for any $p \in
\left[ {2,2 + \frac{{2Q}} {{Q + 2}}\varepsilon _0 } \right)$,
$\Omega ' \subset  \subset \Omega $,
\[
\left\| {Xu} \right\|_{L^p (\Omega ')}  \le c\left( {\left\| g
\right\|_{L_X^{p,\lambda } (\Omega )}  + \left\| {\tilde {\tilde g}}
\right\|_{L_X^{p,\lambda } (\Omega )}  + \left\| {\tilde g}
\right\|_{L_X^{pq_0 ,\lambda q_0 } (\Omega )} } \right).
\]

\textbf{Theorem 1.2} (higher Morrey estimates of gradients for weak
solutions) Under the assumptions in Theorem 1.1, we have that for
any $p \in \left[ {2,2 + \frac{{2Q}} {{Q + 2}}\varepsilon _0 }
\right)$,
\[
Xu \in L_{X,loc}^{p,\lambda } (\Omega ,\mathbb{R}^N ).
\]

\textbf{Theorem 1.3} (H\"{o}lder estimate for weak solutions) Under
the assumptions in Theorem 1.1, it follows that for any $p \in
\left[ {2,2 + \frac{{2Q}} {{Q + 2}}\varepsilon _0 } \right)$, $Q - p
< \lambda < Q$, one has
\[
u \in C_{X,loc}^{0,\kappa } (\Omega ,\mathbb{R}^N ), \quad \kappa  =
1 - \frac{{Q - \lambda }} {p}.
\]

Furthermore, we make the following assumption:

(H3) Let $u \in W_X^{1,2} (\Omega ,\mathbb{R}^N ), g_i (x,u,z)$ and
$f_i^\alpha  (x,u,z)$ satisfy
\[
\left| {g_i (x,u,z)} \right| \le g^i (x) + L\left| z \right|^{\gamma
_0 } ,
\]
\[
\left| {f_i^\alpha  (x,u,z_1 ) - f_i^\alpha  (y,v,z_2 )} \right| \le
L\left( {\left| {g_i^\alpha  (x) - g_i^\alpha  (y)} \right| + \left|
{z_1  - z_2 } \right|} \right),
\]
\[
f_i^\alpha  (x,u,z)z_\alpha ^i  \ge \gamma _1 \left| z \right|^2  -
\left( {g(x)} \right)^2 ,
\]
where $ x,y \in \Omega$, $ u,v \in \mathbb{R}^N$, $z_1,z_2 \in
\mathbb{R}^{qN}$, the selections of $\gamma _0 , L, \gamma _1, g^i,$
$g_i^\alpha$ and $g$, are the same as (H2).

\textbf{Theorem 1.4} (Campanato estimates of gradients for weak
solutions) Let $u \in W_X^{1,2}(\Omega $, $\mathbb{R}^N )$ be a weak
solution to (\ref{eq1}), the coefficients $a_{ij}^{\alpha \beta }$
satisfy (H1), $g_i$ and $f_i^\alpha$ satisfy (H3). Then for any $p
\in \left[ {2,2 + \frac{{2Q}} {{Q + 2}}\varepsilon _0 } \right)$, we
have
\[
Xu \in \mathcal {L}_{X,loc}^{p,\lambda } (\Omega ,\mathbb{R}^N ).
\]

The proof of Theorem 1.1 is based on a priori estimates for weak
solutions to (\ref{eq1}) and the reverse H\"{o}lder inequality in
[13, 27]. In proving Theorem 1.2, several different ways are
attempted and an effective route is the decomposition of (\ref{eq1})
into a nondiagonal homogeneous system and a nondiagonal
nonhomogeneous system. To treat two systems, we discuss regularity
to the homogeneous system corresponding to (\ref{eq1}):
\begin{equation}
\label{eq2}
 - X_\alpha ^ *  (a_{ij}^{\alpha \beta } (x,u)X_\beta  u^j ) = 0.
\end{equation}
With the help of analysis to (\ref{eq2}), we can confirm Theorem 1.2
and Theorem 1.4 under (H2) and (H3), respectively.

Authors in [19, 25] obtained H\"{o}lder regularity for weak
solutions to elliptic systems by employing Liouville theorem.
Differently from this, we prove Theorem 1.3 by combining Morrey
estimates given in Theorem 1.2 and a Morrey lemma in [29].

This paper is organized as follows. In section 2, we introduce
H\"{o}rmander's vector fields, the Carnot-Carath\'{e}odory distance
and some related function spaces, and then recall corresponding
Sobolev-Poincar\'{e} inequality and Morrey Lemma. In section 3, we
prove Theorem 1.1 by choosing appropriate text functions and then
using a priori estimates argument for weak solutions of (\ref{eq1})
and the reverse H\"{o}lder inequality. Section 4 is devoted to the
study of nondiagonal homogeneous degenerate elliptic system
(\ref{eq2}). Through dividing (\ref{eq2}) into two systems which are
a constant coefficients diagonal homogeneous system and a constant
coefficients diagonal nonhomogeneous system, we establish relations
between $L^p$ estimates over balls for gradients of weak solutions
to (\ref{eq2}), see Theorems 4.1 and 4.2. In section 5, we first
divide (\ref{eq1}) into a nondiagonal homogeneous system (5.1) and a
nondiagonal nonhomogeneous system (5.2), and then prove Theorem 1.2
by applying conclusions in sections 3 and 4, and the iteration
lemma. The proof of Theorem 1.3 is given by using Theorem 1.2 and
the known Morrey lemma. After deducing a priori estimates for weak
solutions to (5.2), we finally complete the proof of Theorem 1.4.


\section{Preliminaries}

\label{2}

For every multi-index $\beta  = (\beta _1 ,\beta _2 , \ldots ,\beta
_d )$($1 \le \beta _i  \le q , i = 1, \ldots ,d, \left| \beta
\right| = d$), we call that $d$ is the length of the commutator
$X_\beta   = \left[ {X_{\beta _d } ,\left[ {X_{\beta _{d - 1} } ,
\ldots \left[ {X_{\beta _2 } , \ldots X_{\beta _1 } } \right]}
\right]} \right]$.

\textbf{Definition 2.1} Let $ X_1 , \ldots ,X_q$ be smooth vector
fields. If $\left\{ {X_\beta  \left( {x_0 } \right)}
\right\}_{\left| \beta  \right| \le s}$ spans $\mathbb{R}^n$ at
every $x_0  \in \Omega  \subset \mathbb{R}^n$, then we say that the
system $X = (X_1 , \ldots ,X_q )$ satisfies H\"{o}rmander's
condition of step $s$.

By [26], we can assume that H\"{o}rmander's vector fields $X_1 ,
\ldots ,X_q $ are free up to the order $s$.

\textbf{Definition 2.2} (Carnot-Carath\'{e}odory distance) Let
$\Omega$ be a bounded domain in $\mathbb{R}^n$. An absolutely
continuous curve $\gamma :[0,T] \to \Omega $ is called a sub-unit
curve with respect to the system $X = (X_1 , \ldots,X_q )$, if
$\gamma '(t)$ exists for a.e. $t \in [0,T]$ and satisfies
\[ < \gamma '(t),\xi { > ^2} \le \sum\limits_{j = 1}^q { < {X_j}(\gamma (t)),\xi { > ^2}}
,  \quad \mbox{for any } \xi  \in \mathbb{R}^n. \] We denote the
length of this curve by $l_S \left( \gamma \right) = T$. Given any
$x,y \in \Omega $, let $\Phi (x,y)$ be the collection of all
sub-unit curves connecting $x$ and $y$, and define the
Carnot--Carath\'{e}odory distance induced by $X$ by
\[{d_{\rm{X}}}(x,y) = \inf \{ {l_S}(\gamma ):  \gamma  \in \Phi (x,y)\} .\]

With this distance, we denote a metric ball of radius $R$ centered
at $x_0 $ by
\[{B_R}(x) = B(x,R) = \{ y \in \Omega :d(x,y) < R\} .\]
If one does not need to consider the center of ball, then we also
write $B_R $ instead of $B(x,R)$.

It is well known that the doubling property for metric balls (see
[23]) holds true, i.e., there exist positive constants $c_D $ and
$R_D $, such that for any $x_0 \in \Omega $, $0 < 2R < R_D $,
$B_{2R} \subset \Omega $,
\[ \left| {B(x_0 ,2R)} \right| \le
c_D \left| {B(x_0 ,R)} \right|.
\]
So ${B_R}(x)$ is a homogeneous space ([13]). Furthermore, it follows
that for any $R \le R_D $ and $t \in (0,1)$,
\[\left| {{B_{tR}}} \right| \ge c_D^{ - 1}{t^Q}\left| {{B_R}} \right|.\]
The number $Q = \log _2 c_D $ is called a locally homogeneous
dimension relative to $\Omega $. We can assume by [23] that there
exist two positive constants $c_1 $ and $c_2 $, such that
\begin{equation}
\label{eq3} c_1 R^Q \le \left| {B_R } \right| \le c_2 R^Q.
\end{equation}

\textbf{Definition 2.3} (Sobolev space) Let $1 \le p \le  + \infty,$
$k$ be a positive integer. If $u \in {L^p}(\Omega ,{\mathbb{R}^N})$
satisfies
\[{\left\| u \right\|_{W_X^{k,p}(\Omega ,{\mathbb{R}^N})}} \equiv {\left\| u \right\|_{{L^p}(\Omega ,{\mathbb{R}^N})}}
+ \sum\limits_{h = 1}^k {\sum\limits_{{j_h} = 1}^q {{{\left\|
{{X_{{j_1}}}{X_{{j_2}}} \ldots {X_{{j_h}}}u} \right\|}_{{L^p}(\Omega
,{\mathbb{R}^N})}}} }  <  + \infty, \] then we say that $u$ belongs
to the Sobolev space $W_X^{k,p}(\Omega ,{\mathbb{R}^N})$.

\textbf{Remark: } The space $W_{X,0}^{k,p}(\Omega ,{\mathbb{R}^N})$
is the closure of $C_0^\infty (\Omega ,{\mathbb{R}^N})$ in
$W_X^{k,p}(\Omega ,{\mathbb{R}^N})$ with respect to the norm
${\left\| u \right\|_{W_X^{k,p}(\Omega ,{\mathbb{R}^N})}}$.

Denote by $d_0$ the diameter of $\Omega $.

\textbf{Definition 2.4} (Morrey space) Let $p \ge 1, \lambda  \ge 0,
u \in L_{loc}^p(\Omega ,{\mathbb{R}^N})$, if
\[{\left\| u
\right\|_{L_X^{p,\lambda }(\Omega ,{\mathbb{R}^N})}} \equiv \mathop
{\sup }\limits_{{x_0} \in \Omega ,0 < R < {d_0}} {\left(
{\frac{1}{{{R^\lambda }}}\int_{\Omega  \cap B({x_0},R)} {{{\left|
{u(x)} \right|}^p}dx} } \right)^{\frac{1}{p}}} <  + \infty ,\] then
$u$ is said to belong to the Morrey space $L_X^{p,\lambda }(\Omega
,{\mathbb{R}^N})$.

\textbf{Definition 2.5} (Campanato space) Let $p \ge 1, \lambda  \ge
0, u \in L_{loc}^p(\Omega ,{\mathbb{R}^N})$, if
\[{\left\| u \right\|_{{\mathcal {L}}_X^{p,\lambda }(\Omega ,{\mathbb{R}^N})}}
\equiv \mathop {\sup }\limits_{{x_0} \in \Omega ,0 < R < {d_0}}
{\left( {\frac{1}{{{R^\lambda }}}\int_{\Omega  \cap B({x_0},R)}
{{{\left| {u(x) - {u_{{B_R}}}} \right|}^p}dx} }
\right)^{\frac{1}{p}}} <  + \infty ,\] where ${u_{{B_R}}} =
\frac{1}{{\left| {B({x_0},R)} \right|}}\int_{B({x_0},R)} {u(x)} dx$,
then we say that $u$ is in the Campanato space $\mathcal
{L}_X^{p,\lambda }(\Omega ,{\mathbb{R}^N})$.

\textbf{Definition 2.6} (H\"{o}lder space) Let $\kappa  \in (0,1].$
The H\"{o}lder space $C_X^{0,\kappa }(\bar \Omega ,{\mathbb{R}^N})$
is the set of functions satisfying
\[{\left\| u \right\|_{C_X^{0,\kappa }(\bar \Omega ,
{\mathbb{R}^N})}} \equiv \mathop {\sup }\limits_\Omega  \left| u
\right| + \mathop {\sup }\limits_{\bar \Omega } \frac{{\left| {u(x)
- u(y)} \right|}}{{{{[d(x,y)]}^\kappa }}} <  + \infty. \]

\textbf{Definition 2.7} ($BMO_X$ and $VMO_X$ spaces) Let $u \in
L_{loc}^1(\Omega ,{\mathbb{R}^N})$. If
\[{\left\| u \right\|_{{BMO_X}(\Omega ,{\mathbb{R}^N})}}
 \equiv \mathop {\sup }\limits_{{x_0} \in \Omega ,0 < R < {d_0}} \frac{1}
 {{\left| {\Omega  \cap B({x_0},R)} \right|}}\int_{\Omega  \cap B({x_0},R)}
 {\left| {u(x) - {u_{{B_R}}}} \right|} dx <  + \infty, \]
then we say that $u \in BMO_X(\Omega ,{\mathbb{R}^N})$(Bounded Mean
Oscillation). If $u \in BMO_X(\Omega ,{\mathbb{R}^N})$ and
\[{\eta _R}(u) = \mathop {\sup }\limits_{{x_0} \in \Omega ,0 < \rho  < R} \frac{1}
{{\left| {\Omega  \cap B({x_0},\rho )} \right|}}\int_{\Omega  \cap
B({x_0},\rho )} {\left| {u(x) - {u_{{B_\rho }}}} \right|} dx \to 0,
\quad R \to 0, \] then we say that $u \in VMO_X(\Omega
,{\mathbb{R}^N})$(Vanishing Mean Oscillation).

\textbf{Lemma 2.8 }(see [16]) Let $H(\rho )$ be a nonnegative
increasing function, and for any $0 < \rho  < R \le {R_0} =
dist({x_0},\partial \Omega )$,
\[H(\rho ) \le A\left[ {{{\left( {\frac{\rho }{R}} \right)}^a} + \varepsilon } \right]H(R) + B{R^b},\]
where $A,a$ and $b$ are nonnegative constants with $a>b$. Then there
exist positive constants ${\varepsilon _1} = {\varepsilon
_1}(A,a,b)$ and $c = c(A,a,b)$, such that for any $\varepsilon <
\varepsilon _1 $ , it follows
\[H(\rho ) \le c\left[ {{{\left( {\frac{\rho }{R}} \right)}^b}H(R) + B{\rho ^b}} \right].\]

\textbf{Lemma 2.9} (Sobolev--Poincar\'{e} inequality, see [12] and
[20]) For any open domain $\Omega '$, $\bar \Omega ' \subset \subset
\Omega $, there exist positive constants $R_0$ and $c$, such that
for any $0 < R \le {R_0}$, ${B_R} \subset \Omega $ and $u \in
{C^\infty }(\overline {{B_R}} )$, it holds
\begin{equation}
\label{eq4} {\left( {\frac{1}{{\left| {{B_R}} \right|}}\int_{{B_R}}
{{{\left| {u - {u_R}} \right|}^{p'}}} dx} \right)^{\frac{1}{{p'}}}}
\le cR{\left( {\frac{1}{{\left| {{B_R}} \right|}} \int_{{B_R}}
{{{\left| {Xu} \right|}^p}} dx} \right)^{\frac{1}{p}}},
\end{equation}
where $1 < p < Q, 1 \le p' \le \frac{{pQ}}{{Q - p}}, {u_R} =
\frac{1}{{\left| {{B_R}} \right|}}\int_{{B_R}} {u(x)} dx, R_0$ and
$c$ depend on $\Omega'$ and $\Omega$. In particular, if $u \in
C_0^{^\infty }(\overline {{B_R}} )$, then
\begin{equation}
\label{eq5} {\left( {\frac{1}{{\left| {{B_R}} \right|}}\int_{{B_R}}
{{{\left| u \right|}^{p'}}} dx}
 \right)^{\frac{1}{{p'}}}} \le cR{\left( {\frac{1}{{\left| {{B_R}} \right|}}
 \int_{{B_R}} {{{\left| {Xu} \right|}^p}} dx}
 \right)^{\frac{1}{p}}}.
\end{equation}

\textbf{Lemma 2.10} (Morrey lemma, see [29]) Let $u \in
W_X^{1,p}(\Omega ,{\mathbb{R}^N})(p>1)$ and for any ${B_R} \subset
\subset \Omega $, there exists a constant $\kappa  \in (0,1)$, such
that
$$\int_{{B_R}} {{{\left| {Xu} \right|}^p}} dx \le c{R^{Q - p +
p\kappa }}.$$ Then $u \in C_X^{0,\kappa }(\Omega ,{\mathbb{R}^N})$.


\section{Proof of Theorem 1.1}

\label{3}

The following result is valid to the homogeneous space.

\textbf{Lemma 3.1 }(reverse H\"{o}lder inequality, see [13, 27]) Let
$\hat g,\hat f \ge 0$ satisfy
$$\hat g \in {L^{\hat q}}(\Omega )(\hat q > 1), \hat {f} \in L^{q'}(\Omega )(q' > \hat{q}).$$ Fix a ball
$B_{R_0 } = B(x_0 ,R_0 )$ and assume that for any $x \in B_{R_0 } $
and $R < \frac{1}{2}dist(x,\partial B_{R_0 } )$, there exist
constants $b > 1$ and $\theta \in [0,1)$, such that
\[\frac{1}{{\left| {{B_R}} \right|}}\int_{{B_R}} {{{\hat g}^{\hat q}}dx}
\le b\left[ {{{\left( {\frac{1}{{\left| {{B_{4R/3}}}
\right|}}\int_{{B_{4R/3}}} {\hat gdx} } \right)}^{\hat q}} +
\frac{1}{{\left| {{B_{4R/3}}} \right|}}\int_{{B_{4R/3}}} {{{\hat
f}^{\hat q}}dx} } \right] + \frac{\theta }{{\left| {{B_{4R/3}}}
\right|}} \int_{{B_{4R/3}}} {{{\hat g}^{\hat q}}dx}. \] Then there
exist constants $\varepsilon_0 > 0$ and $c > 0$ such that for any $r
\in [\hat {q}, \hat {q} + \varepsilon_0 )$, it yields $\hat {g} \in
L_{^{loc}}^r (B_{R_0 } )$. Moreover, we have that for any $B_{2R}
\subset \subset \Omega $,
\[{\left( {\frac{1}{{\left| {{B_R}} \right|}}\int_{{B_R}} {{{\hat g}^r}dx} } \right)^{\frac{1}{r}}} \le c\left[ {{{\left( {\frac{1}{{\left| {{B_{2R}}} \right|}}\int_{{B_{2R}}} {{{\hat g}^{\hat q}}dx} } \right)}^{\frac{1}{{\hat q}}}} + {{\left( {\frac{1}{{\left| {{B_{2R}}} \right|}}\int_{{B_{2R}}} {{{\hat f}^r}dx} } \right)}^{\frac{1}{r}}}} \right],
\]
\noindent where $c$ and $\varepsilon_0 $ are positive constants
depending only on $b, \theta , \hat {q}$, and $q'$.

\textbf{Lemma 3.2 } Let the coefficients $a_{ij}^{\alpha \beta }$ in
(1.1) satisfy (H1), functions $g_i$ and $f_i^\alpha$ satisfy (H2).
If $u \in W_X^{1,2}(\Omega ,{\mathbb{R}^N})$ is a weak solution to
(1.1), then for any $p \in \left[ {2,2 + \frac{{2Q}}{{Q +
2}}{\varepsilon _0}} \right)$, where $\varepsilon_0 $ is in Lemma
3.1, ${B_{2R}} \subset \subset \Omega $, we have
\begin{eqnarray}
&& {\left( {\frac{1}{{\left| {{B_R}} \right|}}\int_{{B_R}} {{{\left|
{Xu} \right|}^p}dx} } \right)^{\frac{1}{p}}} \le c{\left(
{\frac{1}{{\left| {{B_{2R}}} \right|}} \int_{{B_{2R}}} {{{\left|
{Xu} \right|}^2}dx} } \right)^{\frac{1}{2}}} \notag \\
&& \quad \quad \quad + c\left[ {{{\left( {\frac{1}{{\left|
{{B_{2R}}} \right|}}\int_{{B_{2R}}} {\left( {{{\left| g \right|}^p}
+ {{\left| {\tilde {\tilde g}} \right|}^p}} \right)dx} }
\right)}^{\frac{1}{p}}} + R{{\left( {\frac{1}{{\left| {{B_{2R}}}
\right|}}\int_{{B_{2R}}} {{{\left| {\tilde g} \right|}^{p{q_0}}}dx}
} \right)}^{\frac{1}{{p{q_0}}}}}} \right]. \label{3.1}
\end{eqnarray}

\textbf{Proof:} We will prove (3.1) with two steps.

\textbf{Step 1.} Let a cut-off function $\eta  \in C_0^\infty
({B_R})$ satisfy that for any $0 < \rho  < R$,
$$0 \le \eta  \le 1 ( \text{ in } B_R), \quad \eta=1 (\text{ in }  B_{\rho}
), \quad \left| {X\eta } \right| \le \frac{c}{{R - \rho }}.$$
Multiplying both sides of (1.1) by the test function $\varphi  = (u
- {u_{{B_R}}}){\eta ^2}$ and integrating on $B_R$, it yields by (H1)
that
\begin{eqnarray}
&& \quad - \int_{{B_R}} {{A^{\alpha \beta }}{\delta _{ij}}{\eta
^2}{X_\beta }{u^j}{X_\alpha }{u^i}dx}  + \int_{{B_R}} {{\eta
^2}f_i^\alpha {X_\alpha }{u^i}dx} \notag \\
&& = \int_{{B_R}} {B_{ij}^{\alpha \beta }{\eta ^2}{X_\beta
}{u^j}{X_\alpha }{u^i}dx}  + \int_{{B_R}} {2a_{ij}^{\alpha \beta
}(x,u)\eta ({u^i} - {u_{{B_R}}}){X_\beta }{u^j}{X_\alpha }\eta dx} \notag \\
&& \quad  + \int_{{B_R}} {\left( {{g_i}({u^i} - {u_{{B_R}}}){\eta
^2} - 2\eta ({u^i} - {u_{{B_R}}})f_i^\alpha {X_\alpha }\eta }
\right)dx}. \label{3.2}
\end{eqnarray}

It shows by (H2), (2.2), H\"{o}lder's inequality and Young's
inequality that
\begin{eqnarray}
&& \quad \int_{B_R } {\left| {g_i } \right|\left| {u^i  - u_R } \right|\eta ^2 dx}
 \le \int_{B_R } {\left( {g^i  + L\left| {Xu} \right|^{\gamma _0 } } \right)
 \left| {u^i  - u_R } \right|dx}  \notag \\
&& \le \left( {\int_{B_R } {\left| {\tilde g} \right|^{\frac{{2Q}}
{{Q + 2}}} dx} } \right)^{\frac{{Q + 2}} {{2Q}}} \left( {\int_{B_R }
{\left| {u - u_R } \right|^{\frac{{2Q}} {{Q - 2}}} dx} }
\right)^{\frac{{Q - 2}}
{{2Q}}}  \notag \\
&& \quad + L\left( {\int_{B_R } {\left| {Xu} \right|^2 dx} }
\right)^{\frac{{\gamma _0 }} {2}} \left( {\int_{B_R } {\left| {u -
u_R } \right|^{\frac{2} {{2 - \gamma _0 }}} dx} } \right)^{\frac{{2
- \gamma _0 }}
{2}}  \notag \\
&& \le \left( {\int_{B_R } {\left| {\tilde g} \right|^{2q_0 } dx} }
\right)^{\frac{1} {{2q_0 }}} \left( {\int_{B_R } {\left| {Xu}
\right|^2 dx} } \right)^{\frac{1}
{2}}  \notag \\
&& \quad + L\left( {\int_{B_R } {\left| {Xu} \right|^2 dx} }
\right)^{\frac{{\gamma _0 }} {2}} cR\left| {B_R } \right|^{\frac{{2
- \gamma _0 }} {2}} \left( {\frac{1} {{\left| {B_R }
\right|}}\int_{B_R } {\left| {Xu} \right|^2 dx} } \right)^{\frac{1}
{2}}  \notag \\
&& \le c_\varepsilon  \left( {\int_{B_R } {\left| {\tilde g}
\right|^{2q_0 } dx} } \right)^{\frac{1} {{q_0 }}}  + \varepsilon
\int_{B_R } {\left| {Xu} \right|^2 dx}  + cR^{\frac{{Q + 2 - Q\gamma
_0 }} {2}} \left( {\int_{B_R } {\left| {Xu} \right|^2 dx} }
\right)^{\frac{{1 + \gamma _0 }}
{2}}  \notag \\
&& \le c_\varepsilon  \left( {\int_{B_R } {\left| {\tilde g}
\right|^{2q_0 } dx} } \right)^{\frac{1} {{q_0 }}}  + \left(
{\varepsilon  + cR^{\frac{{Q + 2 - Q\gamma _0 }} {2}} \left(
{\int_\Omega  {\left| {Xu} \right|^2 dx} } \right)^{\frac{{\gamma _0
- 1}}
{2}} } \right)\int_{B_R } {\left| {Xu} \right|^2 dx} \notag \\
&& \le {c_\varepsilon }{\left( {\int_{{B_R}} {{{\left| {\tilde g}
\right|}^{2{q_0}}}dx} } \right)^{\frac{1}{{{q_0}}}}} + \left(
{\varepsilon  + c{R^{\frac{{Q + 2 - Q{\gamma _0}}}{2}}}}
\right)\int_{{B_R}} {{{\left| {Xu} \right|}^2}dx} .\label{3.3}
\end{eqnarray}
Also,
\begin{eqnarray}
&& \quad \int_{{B_R}} {\eta \left| {{u^i} - {u_R}} \right|\left| {f_i^\alpha } \right|\left| {{X_\alpha }\eta } \right|dx} \notag \\
&& \le \int_{{B_R}} {\eta \left| {{u^i} - {u_R}} \right|\left(
{g_i^\alpha  + L\left| {{X_\alpha }{u^i}} \right|} \right)\left|
{{X_\alpha }\eta } \right|dx} \notag \\
&& \le {c_\varepsilon }\int_{{B_R}} {{{\left| {u - {u_R}}
\right|}^2}{{\left| {X\eta } \right|}^2}dx}  + {c_\varepsilon
}\int_{{B_R}} {{\eta ^2}{{\left| {\tilde {\tilde g}} \right|}^2}dx}
+ \varepsilon \int_{{B_R}} {{\eta ^2}{{\left| {Xu} \right|}^2}dx} .
\label{3.4}
\end{eqnarray}

Inserting (3.3) and (3.4) into (3.2), and noting (H1) and (H2), we
have
\[\begin{array}{l}
 \quad {\lambda _0}\int_{{B_R}} {{\eta ^2}{{\left| {Xu} \right|}^2}dx}
 + {\gamma _1}\int_{{B_R}} {{\eta ^2}{{\left| {Xu} \right|}^2}dx}
 - \int_{{B_R}} {{\eta ^2}{{\left| g \right|}^2}dx}  \\
  \le \delta {\lambda _0}\int_{{B_R}} {{\eta ^2}{{\left| {Xu} \right|}^2}dx}
  + {c_\varepsilon }\int_{{B_R}} {{{\left| {u - {u_R}} \right|}^2}{{\left| {X\eta } \right|}^2}dx}
  + 2\varepsilon \int_{{B_R}} {{\eta ^2}{{\left| {Xu} \right|}^2}dx} \\
\quad + {c_\varepsilon }\int_{{B_R}} {{\eta ^2}{{\left| {\tilde {\tilde g}} \right|}^2}dx}
+ {c_\varepsilon }{\left( {\int_{{B_R}} {{{\left| {\tilde g} \right|}^{2{q_0}}}dx} } \right)^{\frac{1}{{{q_0}}}}}
+ \left( {\varepsilon  + c{R^{\frac{{Q + 2 - Q{\gamma _0}}}{2}}}} \right)\int_{{B_R}} {{{\left| {Xu} \right|}^2}dx} . \\
 \end{array}\]
By properties on $\eta$,
\[\begin{array}{l}
 \quad \left( {{\lambda _0} + {\gamma _1} - 2\varepsilon  - \delta {\lambda _0}} \right)
 \int_{{B_\rho }} {{{\left| {Xu} \right|}^2}dx}  \\
  \le \frac{{{c_\varepsilon }}}{{{{\left( {R - \rho } \right)}^2}}}
  \int_{{B_R}} {{{\left| {u - {u_R}} \right|}^2}dx}
  + {c_\varepsilon }\int_{{B_R}} {\left( {{{\left| g \right|}^2} + {{\left| {\tilde {\tilde g}} \right|}^2}} \right)dx}  \\
\quad + {c_\varepsilon }{\left( {\int_{{B_R}} {{{\left| {\tilde g} \right|}^{2{q_0}}}dx} } \right)^{\frac{1}{{{q_0}}}}} + \left( {\varepsilon  + c{R^{\frac{{Q + 2 - Q{\gamma _0}}}{2}}}} \right)\int_{{B_R}} {{{\left| {Xu} \right|}^2}dx} . \\
 \end{array}\]
Letting $\rho  = \frac{3}{4}R$ and applying (2.2), it obtains
\[\begin{array}{l}
 \quad \frac{{{\lambda _0} + {\gamma _1} - 2\varepsilon  - \delta {\lambda _0}}}
 {{\left| {{B_{3R/4}}} \right|}}\int_{{B_{3R/4}}} {{{\left| {Xu} \right|}^2}dx}  \\
  \le {c_\varepsilon }{\left( {\frac{1}{{\left| {{B_R}} \right|}}
  \int_{{B_R}} {{{\left| {Xu} \right|}^{\frac{{2Q}}{{Q + 2}}}}dx} } \right)^{\frac{{Q + 2}}{Q}}}
  + \frac{{{c_\varepsilon }}}{{\left| {{B_R}} \right|}}\left[ {\int_{{B_R}} {\left( {{{\left| g \right|}^2}
  + {{\left| {\tilde {\tilde g}} \right|}^2}} \right)dx}  + {{\left( {\int_{{B_R}} {{{\left| {\tilde g} \right|}^{2{q_0}}}dx} } \right)}^{\frac{1}{{{q_0}}}}}} \right] \\
 \quad  + \left( {\varepsilon  + c{R^{\frac{{Q + 2 - Q{\gamma _0}}}{2}}}} \right)\frac{1}{{\left| {{B_R}} \right|}}\int_{{B_R}} {{{\left| {Xu} \right|}^2}dx} . \\
 \end{array}\]
Because of $0 < \delta  < 1$, we can choose $\varepsilon$ and $R$
small enough such that ${\lambda _0} + {\gamma _1} - 2\varepsilon -
\delta {\lambda _0} > 0$, and ${\theta _1} = \frac{{\varepsilon  +
c{R^{\frac{{Q + 2 - Q{\gamma _0}}}{2}}}}}{{{\lambda _0} + {\gamma
_1} - 2\varepsilon  - \delta {\lambda _0}}} \in \left( {0,1}
\right)$, so
\begin{eqnarray}
&& \frac{1}{{\left| {{B_{3R/4}}} \right|}}\int_{{B_{3R/4}}}
{{{\left| {Xu} \right|}^2}dx}  \le c{\left( {\frac{1}{{\left|
{{B_R}} \right|}}\int_{{B_R}} {{{\left| {Xu}
\right|}^{\frac{{2Q}}{{Q + 2}}}}dx} } \right)^{\frac{{Q + 2}}{Q}}} +
\frac{{{\theta _1}}}{{\left| {{B_R}} \right|}}\int_{{B_R}} {{{\left|
{Xu} \right|}^2}dx} \notag \\
&& \quad \quad \quad \quad \quad \quad \quad \quad \quad \quad +
\frac{c}{{\left| {{B_R}} \right|}}\left[ {\int_{{B_R}} {\left(
{{{\left| g \right|}^2} + {{\left| {\tilde {\tilde g}} \right|}^2}}
\right)dx}  + {{\left( {\int_{{B_R}} {{{\left| {\tilde g}
\right|}^{2{q_0}}}dx} } \right)}^{\frac{1}{{{q_0}}}}}} \right].
\label{3.5}
\end{eqnarray}

\textbf{Step 2.} Setting $$\hat q = \frac{{Q + 2}}{Q} =
\frac{1}{{{q_0}}}, \hat g = {\left| {Xu} \right|^{\frac{{2Q}}{{Q +
2}}}} = {\left| {Xu} \right|^{2{q_0}}}$$ and $$\hat f = {\left(
{\int_{{B_R}} {{{\left| {\tilde g} \right|}^{2{q_0}}}dx} }
\right)^{1 - {q_0}}}{\left| {\tilde g} \right|^{2{q_0}^2}}+ {\left(
{{{\left| g \right|}^2} + {{\left| {\tilde {\tilde g}} \right|}^2}}
\right)^{{q_0}}},$$ (3.5) can be written as
\[\frac{1}{{\left| {{B_{3R/4}}} \right|}}\int_{{B_{3R/4}}} {{{\hat g}^{\hat q}}dx}
 \le c{\left( {\frac{1}{{\left| {{B_R}} \right|}}\int_{{B_R}} {\hat gdx} } \right)^{\hat q}}
 + \frac{{{\theta _1}}}{{\left| {{B_R}} \right|}}\int_{{B_R}} {{{\hat g}^{\hat q}}dx}
  + \frac{c}{{\left| {{B_R}} \right|}}\int_{{B_R}} {{{\hat f}^{\hat q}}dx} .\]
By Lemma 3.1, we have $\hat g \in L_{loc}^r$, for any $r \in \left[
{\hat q, \hat q + {\varepsilon _0}} \right)$ and ${B_{2R}} \subset
\subset \Omega,$
\[{\left( {\frac{1}{{\left| {{B_R}} \right|}}\int_{{B_R}} {{{\hat g}^r}dx} } \right)^{\frac{1}{r}}}
\le c\left[ {{{\left( {\frac{1}{{\left| {{B_{2R}}} \right|}}
\int_{{B_{2R}}} {{{\hat g}^{\hat q}}dx} } \right)}^{\frac{1}{{\hat
q}}}} + {{\left( {\frac{1}{{\left| {{B_{2R}}}
\right|}}\int_{{B_{2R}}} {{{\hat f}^r}dx} } \right)}^{\frac{1}{r}}}}
\right],\] namely,
\begin{eqnarray}
&& \frac{1}{{\left| {{B_R}} \right|}}\int_{{B_R}} {{{\left| {Xu}
\right|}^{2{q_0}r}}dx}  \le c{\left( {\frac{1}{{\left| {{B_{2R}}}
\right|}}\int_{{B_{2R}}} {{{\left| {Xu} \right|}^2}dx} }
\right)^{{q_0}r}}  \notag \\
&& + \frac{c}{{\left| {{B_{2R}}} \right|}}{\left( {\int_{{B_R}}
{{{\left| {\tilde g} \right|}^{2{q_0}}}dx} } \right)^{\left( {1 -
{q_0}} \right)r}}\int_{{B_{2R}}} {{{\left| {\tilde g}
\right|}^{2{q_0}^2r}}dx}+ \frac{c}{{\left| {{B_{2R}}}
\right|}}\int_{{B_{2R}}} {{{\left( {{{\left| g \right|}^2} +
{{\left| {\tilde {\tilde g}} \right|}^2}} \right)}^{{q_0}r}}dx}.
\label{3.6}
\end{eqnarray}

Denote $p = 2{q_0}r$, then $p \in \left[ {2,2 + \frac{{2Q}}{{Q +
2}}{\varepsilon _0}} \right)$ and
\[
\begin{gathered}
 \quad  \frac{1}
{{\left| {B_R } \right|}}\int_{B_R } {\left| {Xu} \right|^p dx}  \hfill \\
   \le c\left( {\frac{1}
{{\left| {B_{2R} } \right|}}\int_{B_{2R} } {\left| {Xu} \right|^2
dx} } \right)^{\frac{p} {2}}  + \frac{c} {{\left| {B_{2R} }
\right|}}\left( {\int_{B_R } {\left| {\tilde g} \right|^{2q_0 } dx}
} \right)^{\frac{{1 - q_0 }}
{{2q_0 }}p} \int_{B_{2R} } {\left| {\tilde g} \right|^{pq_0 } dx}  \hfill \\
  {\kern 1pt} {\kern 1pt} {\kern 1pt} {\kern 1pt} {\kern 1pt} {\kern 1pt} {\kern 1pt}  + \frac{c}
{{\left| {B_{2R} } \right|}}\int_{B_{2R} } {\left( {\left| g
\right|^2  + \left| {\tilde {\tilde g}} \right|^2 }
\right)^{\frac{p}
{2}} dx}  \hfill \\
   \le c\left( {\frac{1}
{{\left| {B_{2R} } \right|}}\int_{B_{2R} } {\left| {Xu} \right|^2
dx} } \right)^{\frac{p} {2}}  + \frac{c} {{\left| {B_{2R} }
\right|}}\left| {B_{2R} } \right|^{\frac{{(p - 2)(1 - q_0 )}} {{2q_0
}}} \left( {\int_{B_{2R} } {\left| {\tilde g} \right|^{pq_0 } dx} }
\right)^{\frac{{1 - q_0 }}
{{q_0 }}} \int_{B_{2R} } {\left| {\tilde g} \right|^{pq_0 } dx}  \hfill \\
\quad  + \frac{c} {{\left| {B_{2R} } \right|}}\int_{B_{2R} } {\left( {\left| g \right|^p  + \left| {\tilde {\tilde g}} \right|^p } \right)dx}  \hfill \\
   \le c\left[ {\left( {\frac{1}
{{\left| {B_{2R} } \right|}}\int_{B_{2R} } {\left| {Xu} \right|^2
dx} } \right)^{\frac{p} {2}}  + R^p \left( {\frac{1} {{\left|
{B_{2R} } \right|}}\int_{B_{2R} } {\left| {\tilde g} \right|^{pq_0 }
dx} } \right)^{\frac{1}
{{q_0 }}} } \right] \hfill \\
\quad + \frac{c}
{{\left| {B_{2R} } \right|}}\int_{B_{2R} } {\left( {\left| g \right|^p  + \left| {\tilde {\tilde g}} \right|^p } \right)dx} . \hfill \\
\end{gathered}
\]
Hence (3.1) is proved.

\textbf{Corollary 3.3} Let $u \in W_X^{1,2}(\Omega ,{\mathbb{R}^N})$
be a weak solution to the homogeneous degenerate elliptic system
(1.2). Then for any $p \in \left[ {2,2 + \frac{{2Q}}{{Q +
2}}{\varepsilon _0}} \right)$ and ${B_{2R}} \subset \subset \Omega
$, it follows
\begin{eqnarray}
&& {\left( {\frac{1}{{\left| {{B_R}} \right|}}\int_{{B_R}} {{{\left|
{Xu} \right|}^p}dx} } \right)^{\frac{1}{p}}} \le c{\left(
{\frac{1}{{\left| {{B_{2R}}} \right|}}\int_{{B_{2R}}} {{{\left| {Xu}
\right|}^2}dx} } \right)^{\frac{1}{2}}}.\label{3.7}
\end{eqnarray}

\textbf{Proof of theorem 1.1: }Multiplying both sides of (1.1) by
the test function $u - {u_{{B_{2R}}}}$ and integrating on $B_{2R}$,
we have
\[ - \int_{{B_{2R}}} {a_{ij}^{\alpha \beta }(x,u){X_\beta }{u^j}{X_\alpha }({u^i} - {u_{{B_{2R}}}})dx}  = \int_{{B_{2R}}} {\left( {{g_i}({u^i} - {u_{{B_{2R}}}}) - f_i^\alpha {X_\alpha }({u^i} - {u_{{B_{2R}}}})} \right)dx} \]
or
\[\begin{array}{l}
  \quad - \int_{{B_{2R}}} {{A^{\alpha \beta }}{\delta _{ij}}{X_\beta }{u^j}{X_\alpha }{u^i}dx}  + \int_{{B_{2R}}} {f_i^\alpha {X_\alpha }{u^i}dx}  \\
  = \int_{{B_R}} {B_{ij}^{\alpha \beta }{X_\beta }{u^j}{X_\alpha }{u^i}dx}  + \int_{{B_{2R}}} {{g_i}\left( {{u^i} - {u_{{B_{2R}}}}} \right)dx} . \\
 \end{array}\]
By (H1), (H2) and (3.3), it gives
\[\begin{array}{l}
\quad {\lambda _0}\int_{{B_{2R}}} {{{\left| {Xu} \right|}^2}dx}  + {\gamma _1}\int_{{B_{2R}}} {{{\left| {Xu} \right|}^2}dx}  - \int_{{B_{2R}}} {{{\left| g \right|}^2}dx}  \\
  \le \delta {\lambda _0}\int_{{B_R}} {{{\left| {Xu} \right|}^2}dx}  + {\kern 1pt} c\int_{{B_{2R}}} {\left| {{g_i}} \right|\left| {{u^i} - {u_{2R}}} \right|dx}  \\
  \le {c_\varepsilon }{\left( {\int_{{B_{2R}}} {{{\left| {\tilde g} \right|}^{2{q_0}}}dx} } \right)^{\frac{1}{{{q_0}}}}} + \left( {\delta {\lambda _0} + \varepsilon  + c{R^{\frac{{Q + 2 - Q{\gamma _0}}}{2}}}} \right)\int_{{B_{2R}}} {{{\left| {Xu} \right|}^2}dx}  \\
  \buildrel \Delta \over = {c_\varepsilon }{\left( {\int_{{B_{2R}}} {{{\left| {\tilde g} \right|}^{2{q_0}}}dx} } \right)^{\frac{1}{{{q_0}}}}} + {\theta _2}\int_{{B_{2R}}} {{{\left| {Xu} \right|}^2}dx} , \\
 \end{array}\]
where ${\theta _2} = \delta {\lambda _0} + \varepsilon  +
c{R^{\frac{{Q + 2 - Q{\gamma _0}}}{2}}}$. Then
\[\left( {{\lambda _0} + {\gamma _1} - {\theta _2}} \right)\int_{{B_{2R}}} {{{\left| {Xu} \right|}^2}dx}  \le {c_\varepsilon }{\left( {\int_{{B_{2R}}} {{{\left| {\tilde g} \right|}^{2{q_0}}}dx} } \right)^{\frac{1}{{{q_0}}}}} + \int_{{B_{2R}}} {{{\left| g \right|}^2}dx} .\]
Since ${\gamma _0} \in \left[ {1,\frac{{Q + 2}}{Q}} \right)$, $0 <
\delta  < 1$, we can choose $\varepsilon $ and $R$ small enough such
that ${\lambda _0} + {\gamma _1} - {\theta _2} > 0$, and derive
\begin{eqnarray}
&& \int_{{B_{2R}}} {{{\left| {Xu} \right|}^2}dx}  \le c{\left(
{\int_{{B_{2R}}} {{{\left| {\tilde g} \right|}^{2{q_0}}}dx} }
\right)^{\frac{1}{{{q_0}}}}} + c\int_{{B_{2R}}} {{{\left| g
\right|}^2}dx}. \label{3.8}
\end{eqnarray}
It shows from Lemma 3.2 that
\[\begin{array}{l}
 \quad {\left( {\frac{1}{{\left| {{B_R}} \right|}}
 \int_{{B_R}} {{{\left| {Xu} \right|}^p}dx} } \right)^{\frac{1}{p}}} \\
  \le c{\left( {\frac{1}{{\left| {{B_{2R}}} \right|}}
  {{\left( {\int_{{B_{2R}}} {{{\left| {\tilde g} \right|}^{2{q_0}}}dx} } \right)}^{\frac{1}{{{q_0}}}}}
  + \frac{1}{{\left| {{B_{2R}}} \right|}}\int_{{B_{2R}}} {{{\left| g \right|}^2}dx} } \right)^{\frac{1}{2}}} \\
\quad + c{\left( {\frac{1}{{\left| {{B_{2R}}} \right|}}\int_{{B_{2R}}}
 {\left( {{{\left| g \right|}^p} + {{\left| {\tilde {\tilde g}} \right|}^p}} \right)dx} } \right)^{\frac{1}{p}}}
  + R{\left( {\frac{1}{{\left| {{B_{2R}}} \right|}}\int_{{B_{2R}}}
  {{{\left| {\tilde g} \right|}^{p{q_0}}}dx} } \right)^{\frac{1}{{p{q_0}}}}} \\
  \le c\left[ {R{{\left( {\frac{1}{{\left| {{B_{2R}}} \right|}}
  \int_{{B_{2R}}} {{{\left| {\tilde g} \right|}^{p{q_0}}}dx} } \right)}^{\frac{1}{{p{q_0}}}}}
   + {{\left( {\frac{1}{{\left| {{B_{2R}}} \right|}}\int_{{B_{2R}}} {{{\left| g \right|}^p}dx} } \right)}^{\frac{1}{p}}}} \right] \\
\quad + c\left[ {{{\left( {\frac{1}{{\left| {{B_{2R}}} \right|}}\int_{{B_{2R}}} {\left( {{{\left| g \right|}^p}
 + {{\left| {\tilde {\tilde g}} \right|}^p}} \right)dx} } \right)}^{\frac{1}{p}}}
 + R{{\left( {\frac{1}{{\left| {{B_{2R}}} \right|}}\int_{{B_{2R}}}
 {{{\left| {\tilde g} \right|}^{p{q_0}}}dx} } \right)}^{\frac{1}{{p{q_0}}}}}} \right] \\
  \le c\left[ {{{\left( {\frac{1}{{\left| {{B_{2R}}} \right|}}\int_{{B_{2R}}}
  {\left( {{{\left| g \right|}^p} + {{\left| {\tilde {\tilde g}} \right|}^p}} \right)dx} }
   \right)}^{\frac{1}{p}}} + R{{\left( {\frac{1}{{\left| {{B_{2R}}} \right|}}
   \int_{{B_{2R}}} {{{\left| {\tilde g} \right|}^{p{q_0}}}dx} } \right)}^{\frac{1}{{p{q_0}}}}}} \right]. \\
 \end{array}\]
So we conclude by (2.1) that
\begin{eqnarray}
&& \quad \int_{{B_R}} {{{\left| {Xu} \right|}^p}dx}\notag \\
&& \le c\int_{{B_{2R}}} {\left( {{{\left| g \right|}^p}
 + {{\left| {\tilde {\tilde g}} \right|}^p}} \right)dx}
 + c{R^p}\left| {{B_R}} \right|{\left( {\frac{1}{{\left| {{B_{2R}}} \right|}}
 \int_{{B_{2R}}} {{{\left| {\tilde g} \right|}^{p{q_0}}}dx} } \right)^{\frac{1}{{{q_0}}}}} \notag \\
 && \le c\int_{{B_{2R}}} {\left( {{{\left| g \right|}^p} + {{\left| {\tilde {\tilde g}} \right|}^p}} \right)dx}
  + c{R^{p - 2}}{\left( {\int_{{B_{2R}}} {{{\left| {\tilde g} \right|}^{p{q_0}}}dx} } \right)^{\frac{1}{{{q_0}}}}} \notag \\
&& \le c{R^\lambda }\left( {\left\| g \right\|_{L_X^{p,\lambda }}^p
+ \left\| {\tilde {\tilde g}} \right\|_{L_X^{p,\lambda }}^p} \right)
+ c{R^{p + \lambda  - 2}}\left\| {\tilde g} \right\|_{L_X^{p{q_0},\lambda {q_0}}}^p \notag \\
&& \le c{R^\lambda }\left( {\left\| g \right\|_{L_X^{p,\lambda }}^p
+ \left\| {\tilde {\tilde g}} \right\|_{L_X^{p,\lambda }}^p +
\left\| {\tilde g} \right\|_{L_X^{p{q_0},\lambda {q_0}}}^p} \right).
\label{3.9}
\end{eqnarray}
It attains the assertion.

\textbf{Corollary 3.4 } If (H2) in Theorem 1.1 is replaced by (H3),
then the result of Theorem 1.1 still holds.

\textbf{Proof:} Since $u \in W_X^{1,2}(\Omega ,{\mathbb{R}^N})$  is
a weak solution to (1.1), we see that $u$ is also a weak solution to
the following system
\[ - X_\alpha ^ * (a_{ij}^{\alpha \beta }(x,u){X_\beta }{u^j})
= {g_i}(x,u,Xu) - X_\alpha ^ * \left( {f_i^\alpha (x,u,Xu) -
{{\left( {f_i^\alpha } \right)}_{{B_R}}}} \right).\] As in the proof
of Lemma 3.2, it follows
\begin{eqnarray}
&& \quad  - \int_{{B_R}} {{A^{\alpha \beta }}{\delta _{ij}}{\eta
^2}{X_\beta }{u^j}{X_\alpha }{u^i}dx}
+ \int_{{B_R}} {{\eta ^2}f_i^\alpha {X_\alpha }{u^i}dx} \notag \\
&& = \int_{{B_R}} {B_{ij}^{\alpha \beta }{\eta ^2}{X_\beta }{u^j}{X_\alpha }{u^i}dx}
+ \int_{{B_R}} {2a_{ij}^{\alpha \beta }(x,u)\eta ({u^i} - {u_{{B_R}}}){X_\beta }{u^j}{X_\alpha }\eta dx} \notag \\
&& \quad + \int_{{B_R}} {\left( {{g_i}({u^i} - {u_{{B_R}}}){\eta ^2}
- 2\eta ({u^i} - {u_{{B_R}}})\left( {f_i^\alpha  - {{\left(
{f_i^\alpha } \right)}_{{B_R}}}} \right){X_\alpha }\eta } \right)dx}
. \label{3.10}
\end{eqnarray}
Noting (H3), it implies
\begin{eqnarray}
&& \quad \int_{{B_R}} {{{\left| {f_i^\alpha  - {{\left( {f_i^\alpha
} \right)}_{{B_R}}}} \right|}^2}dx}
 \le c\int_{{B_R}} {{{\left| {f_i^\alpha  - \frac{1}{{\left| {{B_R}} \right|}}\int_{{B_R}} {f_i^\alpha dy} } \right|}^2}dx} \notag \\
 && \le c\int_{{B_R}} {{{\left| {\frac{1}{{\left| {{B_R}} \right|}}\int_{{B_R}}
 {\left( {f_i^\alpha \left( {x,u(x),Xu(x)} \right) - f_i^\alpha \left( {y,u(y),Xu(y)} \right)} \right)dy} } \right|}^2}dx} \notag \\
 && \le c\frac{1}{{\left| {{B_R}} \right|}}\int_{{B_R}} {\left( {\int_{{B_R}} {{{\left| {f_i^\alpha \left( {x,u(x),Xu(x)} \right)
 - f_i^\alpha \left( {y,u(y),Xu(y)} \right)} \right|}^2}dy} } \right)dx} \notag \\
 && \le c\frac{1}{{\left| {{B_R}} \right|}}\int_{{B_R}} {\left( {\int_{{B_R}} {{{\left| {g_i^\alpha (x) - g_i^\alpha (y)} \right|}^2}dy}
  + \int_{{B_R}} {{{\left| {Xu(x) - Xu(y)} \right|}^2}dy} } \right)dx} \notag \\
&& \le c\int_{{B_R}} {{{\left| {\tilde{ \tilde g }- {{\left( {\tilde
{\tilde g}} \right)}_{{B_R}}}} \right|}^2}dx}  + c\int_{{B_R}}
{{{\left| {Xu - {{\left( {Xu} \right)}_{{B_R}}}} \right|}^2}dx}
\label{3.11}
\end{eqnarray}
and
\[\begin{array}{l}
\quad \int_{{B_R}} {\eta \left| {{u^i} - {u_{{B_R}}}} \right|\left| {f_i^\alpha  - {{\left( {f_i^\alpha } \right)}_{{B_R}}}} \right|\left| {{X_\alpha }\eta } \right|dx}  \\
  \le {c_\varepsilon }\int_{{B_R}} {{{\left| {u - {u_{{B_R}}}} \right|}^2}{{\left| {X\eta } \right|}^2}dx}
  + \varepsilon \int_{{B_R}} {{\eta ^2}{{\left| {f_i^\alpha  - {{\left( {f_i^\alpha } \right)}_{{B_R}}}} \right|}^2}dx}  \\
  \le {c_\varepsilon }\int_{{B_R}} {{{\left| {u - {u_{{B_R}}}} \right|}^2}{{\left| {X\eta } \right|}^2}dx}
   + \varepsilon \int_{{B_R}} {{{\left| {\tilde {\tilde g }- {{\left( {\tilde{ \tilde g}} \right)}_{{B_R}}}} \right|}^2}dx}  + \varepsilon \int_{{B_R}} {{{\left| {Xu - {{\left( {Xu} \right)}_{{B_R}}}} \right|}^2}dx}  \\
  \le {c_\varepsilon }\int_{{B_R}} {{{\left| {u - {u_{{B_R}}}} \right|}^2}{{\left| {X\eta } \right|}^2}dx}
  + \varepsilon \int_{{B_R}} {{{\left| {\tilde {\tilde g}} \right|}^2}dx}  + \varepsilon \int_{{B_R}} {{{\left| {Xu} \right|}^2}dx} , \\
 \end{array}\]
which indicates that (3.4) still holds. Now we follow the proofs of
Lemma 3.2 and Theorem 1.1 to reach the result.


\section{Homogeneous degenerate elliptic system}

\label{4}

An estimate of gradient of weak solutions to (1.2) is given in
Corollary 3.3. In this section, we continue to study (1.2) and
establish some other useful estimates. The main results in this
section are Theorem 4.1 and Theorem 4.2.

Denote
\[A = {\left( {{A^{\alpha \beta }}(x)} \right)_{{B_R}}} = \frac{1}{{\left| {{B_R}} \right|}}\int_{{B_R}} {{A^{\alpha \beta }}(x)} dx = {S^2} = SS',\]
where $S$ is a positive symmetric matrix, and denote $S({B_R}) =
\left\{ {Sx:x \in {B_R}} \right\}$.

Given a weak solution $u$ to (1.2), let $u = v + w$, here $v \in
W_X^{1,2}\left( {S({B_R}),{\mathbb{R}^N}} \right)$ is a weak
solution to the following constant coefficients diagonal homogeneous
system
\begin{equation}
\label{eq5} \left\{ {\begin{array}{l} - X_\alpha ^ * \left(
{{{\left( {{A^{\alpha \beta }}(x)} \right)}_{{B_R}}}{\delta
_{ij}}{X_\beta }{v^j}} \right) = 0, \\
v - u \in W_{X,0}^{1,2}\left( {S({B_R}),{\mathbb{R}^N}} \right),\\
\end{array}} \right.
\end{equation}
and $w \in W_{X,0}^{1,2}\left( {S({B_R}),{\mathbb{R}^N}} \right)$
satisfies the following constant coefficients diagonal
nonhomogeneous system
\begin{eqnarray}
&&  \quad - X_\alpha ^ * \left( {{{\left( {{A^{\alpha \beta }}(x)}
\right)}_{{B_R}}}{\delta _{ij}}{X_\beta }{w^j}} \right) \notag \\
&&  = - X_\alpha ^ * \left[ {\left( {{{\left( {{A^{\alpha \beta
}}(x)} \right)}_{{B_R}}} - {A^{\alpha \beta }}(x)} \right){\delta
_{ij}}{X_\beta }{u^j}} \right] + X_\alpha ^ * \left( {B_{ij}^{\alpha
\beta }{X_\beta }{u^j}} \right). \label{4.2}
\end{eqnarray}

We start by recalling a lemma in [26].

\textbf{Lemma 4.1 } Let $v_0  \in C^\infty  (\Omega ,\mathbb{R}^N
)$, $B_R  \subset  \subset \Omega$, $k > \frac{Q} {2}$. Then there
exist positive constants $R_0$ and $c$ such that for any $R \le
R_0,$
\begin{eqnarray}
\mathop {\sup }\limits_{x \in B_{R/4} } \left| {v_0 } \right| \le
c\left| {B_R } \right|^{ - \frac{1} {2}} \sum\limits_{\left| I
\right| \leqslant k} {R^{\left| I \right|} \left\| {X_I v_0 }
\right\|_{L^2 (B_R )} } . \label{4.3}
\end{eqnarray}

\textbf{Lemma 4.2} Let $v \in W_X^{1,2} (\Omega ,\mathbb{R}^N )$ be
a weak solution to (4.1). Then $v \in C^\infty  (\Omega )$ and it
follows that for any positive integer $k$ and $S(B_R ) \subset
\subset \Omega$,
\begin{eqnarray}
 \sum\limits_{\left| I \right| \leqslant k}
{\int_{S\left( {B_{R/2^k } } \right)} {\left| {SX_I v} \right|^2 }
dx}  \le \frac{c} {{R^{2k} }}\int_{S(B_R )} {\left| {Sv} \right|^2 }
dx.\label{4.4}
\end{eqnarray}

\textbf{Proof: } Since $A = \left( {A^{\alpha \beta } (x)} \right)_R
= S^2$, it sees that (4.1) can be rewrite as
\[- X_\alpha ^ *  \left( {\delta _{ij} SX_\beta  \left( {Sv^j } \right)} \right) =
 0.\]
By [26], we know that assertions hold.

\textbf{Lemma 4.3} Let $v \in W_X^{1,2} (\Omega ,\mathbb{R}^N )$ be
a weak solution to (4.1). Then for any $0 < \rho  < R $, $S(B_R )
\subset \subset \Omega$,
\begin{eqnarray}
\int_{S(B_\rho  )} {\left| {SXv} \right|^2 } dx \le c\left(
{\frac{\rho } {R}} \right)^Q \int_{S(B_R )} {\left| {SXv} \right|^2
} dx .\label{4.5}
\end{eqnarray}

\textbf{Proof: } Since $u_0 (y) = v(Sy)$ satisfies $ - X_\alpha ^ *
\left( {X_\beta  u_0 ^j } \right) = 0$ in $B_R$, we have by [9] that
\[\int_{B_\rho  } {\left| {Xu_0 } \right|^2 } dy \le c\left(
{\frac{\rho } {R}} \right)^Q \int_{B_R } {\left| {Xu_0 } \right|^2 }
dy.\] By the transformation $x=Sy$, it finishes the proof of (4.5).

\textbf{Theorem 4.1} Let $u \in W_X^{1,2} (\Omega ,\mathbb{R}^N )$
be a weak solution to (1.2) with coefficients $a_{ij}^{\alpha \beta
}$ satisfying (H1). For any $p \in \left[ {2,2 + \frac{{2Q}} {{Q +
2}}\varepsilon _0 } \right)$, $\frac{{(p - 2)Q}} {p} < \mu _1  < Q$,
$0 < \rho  < R$, $S(B_R ) \subset \subset \Omega$, we have
\begin{eqnarray}
\int_{S(B_\rho  )} {\left| {SXu} \right|} ^p dx \le c\left(
{\frac{\rho } {R}} \right)^{\frac{{2Q - p(Q - \mu _1 )}} {2}}
\int_{S(B_R )} {\left| {SXu} \right|} ^p dx. \label{4.6}
\end{eqnarray}

\textbf{Proof:} If $\frac{R} {2} \le \rho  < R$, then the conclusion
is evident. In the sequel it only needs to treat the case $0 < \rho
< \frac{R} {2}$.

First, multiplying both sides in (4.2) by $ w$ and integrating on
$S(B_R )$,
\[
\begin{gathered}
  \quad \int_{S(B_R )} {\left( {A^{\alpha \beta } (x)} \right)_{B_R } \delta _{ij} X_\beta  w^j X_\alpha  w^i dx}  \hfill \\
   = \int_{S(B_R )} {\left( {\left( {A^{\alpha \beta } (x)} \right)_{B_R }  - A^{\alpha \beta } (x)} \right)\delta _{ij} X_\beta  u^j X_\alpha  w^i dx}  \hfill \\
  {\kern 1pt} {\kern 1pt} {\kern 1pt} {\kern 1pt} {\kern 1pt} {\kern 1pt} {\kern 1pt}  - \int_{S(B_R )} {B_{ij}^{\alpha \beta } X_\beta  u^j X_\alpha  w^i dx} . \hfill \\
\end{gathered}
\]
From (H1), Young's inequality and H\"{o}lder's inequality , we have
\begin{eqnarray}
&& \quad \int_{S(B_R )} {\left| {SXw} \right|^2 dx}  \notag \\
&& \le \frac{{c_\varepsilon  }} {{\lambda _0 }}\lambda _0
\int_{S(B_R )} {\left| {\left( {A^{\alpha \beta } (x)} \right)_{B_R
}  - A^{\alpha \beta } (x)} \right|^2 \left| {Xu} \right|^2 dx}  +
\frac{\varepsilon }
{{\lambda _0 }}\lambda _0 \int_{S(B_R )} {\left| {Xw} \right|^2 dx}  \notag \\
&& \quad + \delta \lambda _0 \int_{S(B_R )} {\left| {Xu} \right|\left| {Xw} \right|dx}  \notag \\
&& \le \frac{{c_\varepsilon  }} {{\lambda _0 }}S^2 \int_{S(B_R )}
{\left| {\left( {A^{\alpha \beta } (x)} \right)_{B_R } - \left(
{A^{\alpha \beta } (x)} \right)_{S(B_R )}  + \left( {A^{\alpha \beta
} (x)} \right)_{S(B_R )}
 - A^{\alpha \beta } (x)} \right|^2 \left| {Xu} \right|^2 dx}  \notag \\
 && \quad + \frac{\varepsilon }
{{\lambda _0 }}\int_{S(B_R )} {\left| {SXw} \right|^2 dx}
+ S^2 \delta \int_{S(B_R )} {\left| {Xu} \right|\left| {Xw} \right|dx}  \notag \\
&& \le \frac{{c_\varepsilon  }} {{\lambda _0 }}\int_{S(B_R )}
{\left| {\left( {A^{\alpha \beta } (x)} \right)_{S(B_R )}  -
A^{\alpha \beta } (x)} \right|^2 \left| {SXu} \right|^2 dx}  +
\frac{{c_\varepsilon  }}
{{\lambda _0 }}\int_{S(B_R )} {\left| {SXu} \right|^2 dx}  \notag \\
 && \quad + \frac{\varepsilon }
{{\lambda _0 }}\int_{S(B_R )} {\left| {SXw} \right|^2 dx}
+ \delta \int_{S(B_R )} {\left| {SXu} \right|\left| {SXw} \right|dx}  \notag \\
&& \le \frac{{c_\varepsilon  }} {{\lambda _0 }}\left( {\int_{S(B_R
)} {\left| {\left( {A^{\alpha \beta } (x)} \right)_{S(B_R )}  -
A^{\alpha \beta } (x)} \right|^{\frac{{2p}} {{p - 2}}} dx} }
\right)^{\frac{{p - 2}} {p}} \left( {\int_{S(B_R )}
{\left| {SXu} \right|^p } dx} \right)^{\frac{2} {p}} \notag \\
&& \quad + \left( {\frac{\varepsilon } {{\lambda _0 }} + \varepsilon
} \right)\int_{S(B_R )} {\left| {SXw} \right|^2 dx}  + \left(
{\frac{{c_\varepsilon  }} {{\lambda _0 }} + c_\varepsilon  }
\right)\int_{S(B_R )} {\left| {SXu} \right|^2 dx} . \label{4.7}
\end{eqnarray}
Noting
\[
\left( {\int_{S(B_R )} {\left| {\left( {A^{\alpha \beta } }
\right)_{S(B_R )}  - A^{\alpha \beta } } \right|^{\frac{{2p}} {{p -
2}}} dx} } \right)^{\frac{{p - 2}} {p}}  \le c\left| {S(B_R )}
\right|^{\frac{{p - 2}} {p}} \left( {\eta _{S(B_R )} \left(
{A^{\alpha \beta } } \right)} \right)^{\frac{{p - 2}} {p}},
\]
it obtains
\[\begin{array}{l}
\quad \int_{S({B_R})} {{{\left| {SXw} \right|}^2}dx}  \\
  \le \frac{{{c_\varepsilon }}}{{{\lambda _0}}}{\left| {S({B_R})} \right|^{\frac{{p - 2}}{p}}}{\left( {{\eta _{S({B_R})}}\left( {{A^{\alpha \beta }}} \right)} \right)^{\frac{{p - 2}}{p}}}{\left( {\int_{S({B_R})} {{{\left| {SXu} \right|}^p}} dx} \right)^{\frac{2}{p}}} + \left( {\frac{\varepsilon }{{{\lambda _0}}} + \varepsilon } \right)\int_{S({B_R})} {{{\left| {SXw} \right|}^2}dx}  \\
 {\kern 1pt} {\kern 1pt} {\kern 1pt} {\kern 1pt} {\kern 1pt} {\kern 1pt} {\kern 1pt} {\kern 1pt} {\kern 1pt}  + \left( {\frac{{{c_\varepsilon }}}{{{\lambda _0}}} + {c_\varepsilon }} \right){\left| {S({B_R})} \right|^{\frac{{p - 2}}{p}}}{\left( {\int_{S({B_R})} {{{\left| {SXu} \right|}^p}} dx} \right)^{\frac{2}{p}}} \\
  \le {c_\varepsilon }{\left| {S({B_R})} \right|^{\frac{{p - 2}}{p}}}\left[ {{{\left( {{\eta _{S({B_R})}}\left( {{A^{\alpha \beta }}} \right)} \right)}^{\frac{{p - 2}}{p}}} + {c_\varepsilon }} \right]{\left( {\int_{S({B_R})} {{{\left| {SXu} \right|}^p}} dx} \right)^{\frac{2}{p}}} \\
 {\kern 1pt} {\kern 1pt} {\kern 1pt} {\kern 1pt} {\kern 1pt} {\kern 1pt} {\kern 1pt} {\kern 1pt}  + \left( {\frac{\varepsilon }{{{\lambda _0}}} + \varepsilon } \right)\int_{S({B_R})} {{{\left| {SXw} \right|}^2}dx} . \\
 \end{array}\]
Choosing $\varepsilon$ small enough such that $1 - \frac{\varepsilon
} {{\lambda _0 }} - \varepsilon  > 0$, it follows
\begin{eqnarray}
\int_{S(B_R )} {\left| {SXw} \right|^2 dx}  \le c\left| {S(B_R )}
\right|\left[ {\left( {\eta _{S(B_R )} \left( {A^{\alpha \beta } }
\right)} \right)^{\frac{{p - 2}} {p}}  + c} \right]\left( {\frac{1}
{{\left| {S(B_R )} \right|}}\int_{S(B_R )} {\left| {SXu} \right|^p }
dx} \right)^{\frac{2} {p}} . \label{4.8}
\end{eqnarray}

Next by $u=v+w$, we have
\begin{eqnarray}
&& \quad \int_{S(B_{2\rho } )} {\left| {SXu} \right|^2 dx} \le
2\int_{S(B_{2\rho } )} {\left| {SXv} \right|^2 dx}
+ 2\int_{S(B_{2\rho } )} {\left| {SXw} \right|^2 dx}  \notag \\
&& \le c\left( {\frac{\rho }
{R}} \right)^Q \int_{S(B_R )} {\left| {SXu} \right|^2 dx}  + c\int_{S(B_R )} {\left| {SXw} \right|^2 dx}  \notag \\
 && \le c\left( {\frac{\rho }
{R}} \right)^Q \left| {S(B_R )} \right|\left( {\frac{1} {{\left|
{S(B_R )} \right|}}\int_{S(B_R )} {\left| {SXu} \right|^p dx} }
\right)^{\frac{2}
{p}}  \notag \\
&&\quad  + c\left| {S(B_R )} \right|\left[ {\left( {\eta _{S(B_R )}
\left( {A^{\alpha \beta } } \right)} \right)^{\frac{{p - 2}} {p}}  +
c} \right]\left( {\frac{1} {{\left| {S(B_R )} \right|}}\int_{S(B_R
)} {\left| {SXu} \right|^p } dx} \right)^{\frac{2}
{p}}  \notag \\
&& \le c\left| {S(B_R )} \right|\left( {\left( {\frac{\rho } {R}}
\right)^Q  + \left( {\eta _{S(B_R )} \left( {A^{\alpha \beta } }
\right)} \right)^{\frac{{p - 2}} {p}}  + c} \right)\left( {\frac{1}
{{\left| {S(B_R )} \right|}}\int_{S(B_R )} {\left| {SXu} \right|^p }
dx} \right)^{\frac{2} {p}} .\label{4.9}
\end{eqnarray}
Using (3.7), we have
\[
\left( {\frac{1} {{\left| {S(B_\rho  )} \right|}}\int_{S(B_\rho  )}
{\left| {SXu} \right|^p dx} } \right)^{\frac{1} {p}}  \le c\left(
{\frac{1} {{\left| {S(B_{2\rho } )} \right|}}\int_{S(B_{2\rho } )}
{\left| {Xu} \right|^2 dx} } \right)^{\frac{1} {2}} .
\]
Inserting (4.9) into the above, it gets
\[
\begin{gathered}
 \quad \int_{S(B_\rho  )} {\left| {SXu} \right|^p dx}  \hfill \\
   \le c\left( {\left( {\frac{\rho }
{R}} \right)^Q  + \left( {\eta _{S(B_R )} \left( {A^{\alpha \beta }
} \right)} \right)^{\frac{{p - 2}} {p}}  + c} \right)^{\frac{p} {2}}
\left( {\frac{{\left| {S(B_\rho  )} \right|}} {{\left| {S(B_R )}
\right|}}} \right)^{\frac{{2 - p}}
{2}} \int_{S(B_R )} {\left| {SXu} \right|^p dx} . \hfill \\
\end{gathered}
\]
Therefore,
\[
\begin{gathered}
 \quad \left( {\left| {S(B_\rho  )} \right|^{\frac{{p - 2}}
{2}} \int_{S(B_\rho  )} {\left| {SXu} \right|^p dx} }
\right)^{\frac{2}
{p}}  \hfill \\
   \le c\left( {\left( {\frac{\rho }
{R}} \right)^Q  + \left( {\eta _{S(B_R )} \left( {A^{\alpha \beta }
} \right)} \right)^{\frac{{p - 2}} {p}}  + c} \right)\left( {\left|
{S(B_R )} \right|^{\frac{{p - 2}} {2}} \int_{S(B_R )} {\left| {SXu}
\right|^p dx} } \right)^{\frac{2}
{p}} . \hfill \\
\end{gathered}
\]

Finally, let $$H(\rho ) = \left( {\left| {S(B_\rho  )}
\right|^{\frac{{p - 2}} {2}} \int_{S(B_\rho  )} {\left| {SXu}
\right|^p dx} } \right)^{\frac{2} {p}},$$ $$H(R) = \left( {\left|
{S(B_R )} \right|^{\frac{{p - 2}} {2}} \int_{S(B_R )} {\left| {SXu}
\right|^p dx} } \right)^{\frac{2} {p}},$$ $$a = Q, \quad B = 0.$$
For any $\mu _1$, $\frac{{(p - 2)Q}} {p} < \mu _1  < Q$, let $ b =
\mu _1$, then $a>b$. Now we apply Lemma 2.8 to reach
\[
\left( {\left| {S(B_\rho  )} \right|^{\frac{{p - 2}} {2}}
\int_{S(B_\rho  )} {\left| {SXu} \right|^p dx} } \right)^{\frac{2}
{p}}  \le c\left( {\frac{\rho } {R}} \right)^{\mu _1 } \left(
{\left| {S(B_R )} \right|^{\frac{{p - 2}} {2}} \int_{S(B_R )}
{\left| {SXu} \right|^p dx} } \right)^{\frac{2} {p}}
\]
and (4.6) is proved.

\textbf{Lemma 4.4 } Let $ w \in W_{X,0}^{1,2} (\Omega ,\mathbb{R}^N
)$ be a weak solution to (4.2), with the coefficients
$a_{ij}^{\alpha \beta }$ satisfying (H1). Then for any $p \in \left[
{2,2 + \frac{{2Q}} {{Q + 2}}\varepsilon _0 } \right)$, $0 < \rho <
R$, $S(B_R ) \subset \subset \Omega$,
\begin{eqnarray}
\int_{S(B_R )} {\left| {SXw} \right|^p } dx \le c\int_{S(B_{2R} )}
{\left| {SXu - \left( {SXu} \right)_{S(B_{2R} )} } \right|^p } dx.
\label{4.10}
\end{eqnarray}

\textbf{Proof: } Clearly, $w$ is also a weak solution to the
following system
\begin{eqnarray}
 - X_\alpha ^ *  \left( {\left( {A^{\alpha \beta } (x)} \right)_{B_R }
 \delta _{ij} X_\beta  w^j } \right) =
- X_\alpha ^ *  \left[ {\left( {A^{\alpha \beta } (x)} \right)_{B_R
} \delta _{ij} \left( {X_\beta  u^j  - \left( {X_\beta  u^j }
\right)_{B_{2R} } } \right)} \right]. \label{4.11}
\end{eqnarray}
Take the cut-off function $\eta  \in C_0^\infty  (B_R )$ as in the
proof of Lemma 3.2. Multiplying both sides of (4.11) by $\varphi  =
(w - w_{B_R } )\eta ^2 $ and integrating on $S(B_R )$, it gets
\[
\begin{gathered}
 \quad \int_{S(B_R )} {\left( {A^{\alpha \beta } (x)} \right)_{B_R } \delta _{ij} X_\beta  w^j X_\alpha  \left( {(w^i  - w_{B_R } )\eta ^2 } \right)dx}  \hfill \\
   = \int_{S(B_R )} {\left( {A^{\alpha \beta } (x)} \right)_{B_R } \delta _{ij} \left( {X_\beta  u^j  - \left( {X_\beta  u^j } \right)_{B_{2R} } } \right)X_\alpha  \left( {(w^i  - w_{B_R } )\eta ^2 } \right)dx} , \hfill \\
\end{gathered}
\]
i.e.,
\[
\begin{gathered}
 \quad \int_{S(B_R )} {\left( {A^{\alpha \beta } (x)} \right)_{B_R } \delta _{ij} \eta ^2 X_\alpha  w^i X_\beta  w^i } dx \hfill \\
   =  - \int_{S(B_R )} {2\left( {A^{\alpha \beta } (x)} \right)_{B_R } \delta _{ij} \eta (w^i  - w_{B_R } )X_\alpha  \eta X_\beta  w^j } dx \hfill \\
  {\kern 1pt} {\kern 1pt} {\kern 1pt} {\kern 1pt} {\kern 1pt} {\kern 1pt} {\kern 1pt} {\kern 1pt}  + \int_{S(B_R )} {\left( {A^{\alpha \beta } (x)} \right)_{B_R } \delta _{ij} \eta ^2 X_\alpha  w^i \left( {X_\beta  u^j  - \left( {X_\beta  u^j } \right)_{B_{2R} } } \right)dx}  \hfill \\
  {\kern 1pt} {\kern 1pt} {\kern 1pt} {\kern 1pt} {\kern 1pt} {\kern 1pt} {\kern 1pt}  + \int_{S(B_R )} {2\left( {A^{\alpha \beta } (x)} \right)_{B_R } \delta _{ij} \eta (w^i  - w_{B_R } )X_\alpha  \eta \left( {X_\beta  u^j  - \left( {X_\beta  u^j } \right)_{B_{2R} } } \right)dx} . \hfill \\
\end{gathered}
\]
It yields by (H1) and Young's inequality that
\[\begin{array}{l}
 \int_{S({B_R})} {{\eta ^2}{{\left| {SXw} \right|}^2}} dx \le {c_\varepsilon }
 \int_{S({B_R})} {{{\left| {Sw - S{w_{{B_R}}}} \right|}^2}{{\left| {X\eta } \right|}^2}} dx \\
 \quad \quad \quad \quad  + {c_\varepsilon }\int_{S({B_R})} {{\eta ^2}{{\left| {SXu - S{{\left( {Xu} \right)}_{{B_{2R}}}}} \right|}^2}dx}  + 2\varepsilon \int_{S({B_R})} {{\eta ^2}{{\left| {SXw} \right|}^2}dx} . \\
 \end{array}\]
From properties of $\eta$ and $Sw_{B_R }  = \left( {Sw}
\right)_{S(B_R )}$, we have
\[\begin{array}{l}
 \int_{S({B_\rho })} {{{\left| {SXw} \right|}^2}} dx \le \frac{{{c_\varepsilon }}}{{{{(R - \rho )}^2}}}
 \int_{S({B_R})} {{{\left| {Sw - {{\left( {Sw} \right)}_{S({B_R})}}} \right|}^2}} dx \\
   \quad \quad  \quad \quad \quad \quad  + {c_\varepsilon }\int_{S({B_R})} {{{\left| {SXu - {{\left( {SXu} \right)}_{S({B_{2R}})}}} \right|}^2}dx}  + 2\varepsilon \int_{S({B_R})} {{{\left| {SXw} \right|}^2}dx} {\kern 1pt} . \\
 \end{array}\]
Letting $ \rho  = \frac{3} {4}R$ and using (2.2), it follows
\begin{eqnarray}
&&\quad \int_{S(B_{3R/4} )} {\left| {SXw} \right|^2 } dx \notag \\
 && \le \frac{{c_\varepsilon  }}
{{R^2 }}\int_{S(B_R )} {\left| {Sw - \left( {Sw} \right)_{S(B_R )} }
\right|^2 } dx + c_\varepsilon  \int_{S(B_R )} {\left| {SXu - \left(
{SXu} \right)_{S(B_{2R} )} } \right|^2 dx} \notag \\
&& \quad + 2\varepsilon \int_{S(B_R )} {\left| {SXw} \right|^2 dx}   \notag \\
&& \le c_\varepsilon  \left| {S(B_R )} \right|\left( {\frac{1}
{{\left| {S(B_R )} \right|}}\int_{S(B_R )} {\left| {SXw}
\right|^{\frac{{2Q}} {{Q + 2}}} } dx} \right)^{\frac{{Q + 2}}
{Q}}   \notag \\
&& \quad  + c_\varepsilon  \int_{S(B_R )} {\left| {SXu - \left(
{SXu} \right)_{S(B_{2R} )} } \right|^2 dx} + 2\varepsilon
\int_{S(B_R )} {\left| {SXw} \right|^2 dx} . \label{4.12}
\end{eqnarray}

Drawing notations $\hat g = \left| {SXw} \right|^{\frac{{2Q}} {{Q +
2}}}$, $\hat q = \frac{{Q + 2}} {Q}$ and $ \hat f = \left| {SXu -
(SXu)_{S(B_{2R} )} } \right|^{\frac{{2Q}} {{Q + 2}}}$, the
inequality (4.12) is of the form
\begin{eqnarray}
&& \quad \frac{1} {{\left| {S(B_{3R/4} )} \right|}}\int_{S(B_{3R/4}
)}{\hat g^{\hat q} dx} \notag \\
&& \le c\left( {\frac{1} {{\left| {S(B_R )} \right|}}\int_{S(B_R )}
{\hat gdx} } \right)^{\hat q}  + \frac{c} {{\left| {S(B_R )}
\right|}}\int_{S(B_R )} {\hat f^{\hat q} dx}  + \frac{{2\varepsilon
}} {{\left| {S(B_R )} \right|}}\int_{S(B_R )} {\hat g^{\hat q} dx} .
\label{4.13}
\end{eqnarray}
We know by choosing $\varepsilon $ small enough such that
$2\varepsilon < 1$ and employing Lemma 3.1 that $\hat g \in
L_{loc}^r$, $r \in \left[ {\hat q,\hat q + \varepsilon _0 }
\right)$, and for any $S(B_{2R} ) \subset  \subset \Omega$,
\begin{eqnarray}
&&{\left( {\frac{1}{{\left| {S({B_R})} \right|}}\int_{S({B_R})}
{{{\left| {SXw} \right|}^{\frac{{2Qr}}{{Q + 2}}}}dx} }
\right)^{\frac{1}{r}}} \le c{\left( {\frac{1}{{\left| {S({B_{2R}})}
\right|}}\int_{S({B_{2R}})} {{{\left| {SXw} \right|}^2}dx} }
\right)^{\frac{Q}{{Q + 2}}}} \notag \\
&& \quad \quad \quad \quad \quad \quad \quad \quad + c{\left(
{\frac{1}{{\left| {S({B_{2R}})} \right|}}\int_{S({B_{2R}})}
{{{\left| {SXu - {{(SXu)}_{S({B_{2R}})}}} \right|}^{\frac{{2Qr}}{{Q
+ 2}}}}dx} } \right)^{\frac{1}{r}}}. \label{4.14}
\end{eqnarray}
Let $ p = \frac{{2Qr}} {{Q + 2}}$, then $p \in \left[ {2,2 +
\frac{{2Q}} {{Q + 2}}\varepsilon _0 } \right)$ and we can rewrite
(4.14) as
\begin{eqnarray}
&& {\left( {\frac{1}{{\left| {S({B_R})} \right|}}\int_{S({B_R})}
{{{\left| {SXw} \right|}^p}dx} } \right)^{\frac{1}{p}}} \le c{\left(
{\frac{1}{{\left| {S({B_{2R}})} \right|}}\int_{S({B_{2R}})}
{{{\left| {SXw} \right|}^2}dx} } \right)^{\frac{1}{2}}} \notag \\
&& \quad \quad \quad \quad \quad \quad \quad \quad + c{\left(
{\frac{1}{{\left| {S({B_{2R}})} \right|}}\int_{S({B_{2R}})}
{{{\left| {SXu - {{(SXu)}_{S({B_{2R}})}}} \right|}^p}dx} }
\right)^{\frac{1}{p}}}. \label{4.15}
\end{eqnarray}

On the other hand, multiplying both sides of (4.11) by $w$ and
integrating on $S(B_{2R} )$,
\[\begin{array}{l}
\quad \int_{S({B_{2R}})} {{{\left( {{A^{\alpha \beta }}(x)} \right)}_{{B_R}}}{\delta _{ij}}{X_\beta }{w^j}{X_\alpha }{w^i}dx}  \\
  = \int_{S({B_{2R}})} {{{\left( {{A^{\alpha \beta }}(x)} \right)}_{{B_R}}}{\delta _{ij}}\left( {{X_\beta }{u^j} - {{\left( {{X_\beta }{u^j}} \right)}_{{B_{2R}}}}} \right){X_\alpha }{w^i}dx} . \\
 \end{array}\]
It follows by (H1) and Young's inequality that
\[
\int_{S(B_{2R} )} {\left| {SXw} \right|^2 dx}  \le c_\varepsilon
\int_{S(B_{2R} )} {\left| {SXu - \left( {SXu} \right)_{S(B_R )} }
\right|^2 dx + \varepsilon \int_{S(B_{2R} )} {\left| {SXw} \right|^2
dx} } .
\]
For $\varepsilon$ small enough, we see
\[
\int_{S(B_{2R} )} {\left| {SXw} \right|^2 dx}  \le c_\varepsilon
\int_{S(B_{2R} )} {\left| {SXu - \left( {SXu} \right)_{S(B_R )} }
\right|^2 dx} .
\]
Putting it into (4.15) implies
\[
\begin{gathered}
 \quad \frac{1} {{\left| {S(B_R )} \right|}}\int_{S(B_R )} {\left| {SXw} \right|^p dx}  \hfill \\
   \le c\left( {\frac{1}
{{\left| {S(B_{2R} )} \right|}}\int_{S(B_{2R} )} {\left| {SXu -
\left( {SXu} \right)_{S(B_R )} } \right|^2 dx} } \right)^{\frac{p}
{2}}  \hfill \\
 \quad + \frac{c} {{\left| {S(B_{2R} )} \right|}}\int_{S(B_{2R} )} {\left| {SXu - (SXu)_{S(B_{2R} )} } \right|^p dx}  \hfill \\
   \le \frac{c} {{\left| {S(B_{2R} )} \right|}}\int_{S(B_{2R} )} {\left| {SXu - (SXu)_{S(B_{2R} )} } \right|^p dx} . \hfill \\
\end{gathered}
\]
The proof of (4.10) is ended.

\textbf{Lemma 4.5 } Let $ v \in W_{X}^{1,2} (\Omega ,\mathbb{R}^N )$
be a weak solution to (4.1). Then for any $p \in \left[ {2,2 +
\frac{{2Q}} {{Q + 2}}\varepsilon _0 } \right)$, $0 < \rho < R$,
$S(B_R ) \subset \subset \Omega$, we have
\begin{eqnarray}
\int_{S(B_\rho  )} {\left| {SXv - \left( {SXv} \right)_{S(B_\rho  )}
} \right|^p } dx \le c\left( {\frac{\rho } {R}} \right)^{Q + p}
\int_{S(B_R )} {\left| {SXv - \left( {SXv} \right)_{S(B_R )} }
\right|^p } dx. \label{4.16}
\end{eqnarray}

\textbf{Proof: } Let $k$ be a fixed integer such that $k > \frac{Q}
{2}$. If $\frac{R} {{2^{k + 2} }} \le \rho  < R$, then the
conclusion is evident. If $\rho  < \frac{R} {{2^{k + 2} }} $, then
$Xv$ and $X^2 v$ are also weak solutions to (4.1), so (4.3) is true
for $X^2 v$ and (4.4) is true for $Xv$. Combining these and noting
(2.1), it shows
\begin{eqnarray}
&& \quad \int_{S(B_\rho  )} {\left| {SX^2 v} \right|^p } dx \le
\left| {S(B_\rho  )} \right|\mathop {\sup }
\limits_{S\left( {B_{R/2^{k + 2} } } \right)} \left| {SX^2 v} \right|^p  \notag \\
&& \le c\left| {S(B_\rho  )} \right|\sum\limits_{\left| I \right|
\le k} {\left| {S\left( {B_{R/2^k } } \right)} \right|^{ - \frac{p}
{2}} R^{p\left| I \right|} \left( {\int_{S\left( {B_{R/2^k } }
\right)} {\left| {SX_I X^2 v} \right|^2 } dx} \right)} ^{\frac{p}
{2}}  \notag \\
 && \le c\left| {S(B_\rho  )} \right|\left| {S(B_R )} \right|^{ - \frac{p}
{2}} \sum\limits_{\left| I \right| \leqslant k} {R^{p\left| I
\right|} \left( {R^{ - 2(\left| I \right| + 1)} \int_{S\left( {B_R }
\right)} {\left| {SXv} \right|^2 } dx} \right)^{\frac{p}
{2}} }  \notag  \\
&& \le c\left| {S(B_\rho  )} \right|\left| {S(B_R )} \right|^{ -
\frac{p} {2}} \sum\limits_{\left| I \right| \leqslant k} {R^{p\left|
I \right|} R^{ - p(\left| I \right| + 1)} \left( {\int_{S(B_R )}
{\left| {SXv} \right|^2 } dx} \right)^{\frac{p}
{2}} }  \notag \\
 && \le c\left| {S(B_\rho  )} \right|\left| {S(B_R )} \right|^{ - \frac{p}
{2}} R^{ - p} \left| {S(B_R )} \right|^{\frac{{p - 2}}
{2}} \int_{S(B_R )} {\left| {SXv} \right|^p } dx \notag \\
 && \le c\frac{{\left| {S(B_\rho  )} \right|}}
{{\left| {S(B_R )} \right|}}R^{ - p} \int_{S(B_R )} {\left| {SXv} \right|^p } dx \notag \\
&& \le c\left( {\frac{\rho } {R}} \right)^Q R^{ - p} \int_{S(B_R )}
{\left| {SXv} \right|^p } dx. \label{4.17}
\end{eqnarray}
Since $ Xv - \left( {Xv} \right)_{B_R }$ is a weak solution to
(4.1), we know that (4.17) is valid for $ Xv - \left( {Xv}
\right)_{B_R }$ and then
\begin{eqnarray}
\int_{S(B_\rho  )} {\left| {SX^2 v} \right|^p } dx \le c\left(
{\frac{\rho } {R}} \right)^Q R^{ - p} \int_{S(B_R )} {\left| {SXv -
\left( {SXv} \right)_{S(B_R )} } \right|^p } dx. \label{4.18}
\end{eqnarray}
Using (2.2) and (4.18), it follows
\[\begin{array}{l}
 \quad \int_{S({B_\rho })} {{{\left| {SXv - {{\left( {SXv} \right)}_{S({B_\rho })}}} \right|}^p}} dx \le c{\rho ^p}\int_{S({B_\rho })} {{{\left| {S{X^2}v} \right|}^p}} dx \\
  \le c{\left( {\frac{\rho }{R}} \right)^{Q + p}}\int_{S({B_R})} {{{\left| {SXv - {{\left( {SXv} \right)}_{S({B_R})}}} \right|}^p}} dx \\
 \end{array}\]
and (4.16) is proved.

\textbf{Theorem 4.2} Let $u \in W_X^{1,2} (\Omega ,\mathbb{R}^N )$
be a weak solution to (1.2), with the coefficients $a_{ij}^{\alpha
\beta }$ satisfying (H1). Then for any $p \in \left[ {2,2 +
\frac{{2Q}} {{Q + 2}}\varepsilon _0 } \right)$, $0 < \mu _2 < Q+p $,
$0 < \rho < R$, $S(B_R ) \subset \subset \Omega$,
\begin{eqnarray}
\int_{S(B_\rho  )} {\left| {SXu - \left( {SXu} \right)_{S(B_\rho )}
} \right|} ^p dx \le c\left( {\frac{\rho } {R}} \right)^{\mu _2 }
\int_{S(B_R )} {\left| {SXu - \left( {SXu} \right)_{S(B_R )} }
\right|} ^p dx. \label{4.19}
\end{eqnarray}

\textbf{Proof: } Noting Lemma 4.4, Lemma 4.5 and $u=v+w$, it leads
to
$$\begin{array}{l}
 \quad {\int_{S({B_\rho })} {\left| {SXu - {{\left( {SXu} \right)}_{S({B_\rho })}}} \right|} ^p}dx \\
  \le c{\int_{S({B_\rho })} {\left| {SXv - {{\left( {SXv} \right)}_{S({B_\rho })}}} \right|} ^p}dx + c{\int_{S({B_\rho })} {\left| {SXw - {{\left( {SXw} \right)}_{S({B_\rho })}}} \right|} ^p}dx \\
  \le c{\left( {\frac{\rho }{R}} \right)^{Q + p}}{\int_{S({B_R})} {\left| {SXv - {{\left( {SXv} \right)}_{S({B_R})}}} \right|} ^p}dx +  + c{\int_{S({B_\rho })} {\left| {SXw} \right|} ^p}dx \\
  \le c{\left( {\frac{\rho }{R}} \right)^{Q + p}}{\int_{S({B_R})} {\left| {SXu - {{\left( {SXu} \right)}_{S({B_R})}}} \right|} ^p}dx + c{\int_{S({B_R})} {\left| {SXw - {{\left( {SXw} \right)}_{S({B_R})}}} \right|} ^p}dx \\
 {\kern 1pt} {\kern 1pt} {\kern 1pt} {\kern 1pt} {\kern 1pt} {\kern 1pt}  + c{\int_{S({B_\rho })} {\left| {SXw} \right|} ^p}dx \\
  \le c{\left( {\frac{\rho }{R}} \right)^{Q + p}}{\int_{S({B_R})} {\left| {SXu - {{\left( {SXu} \right)}_{S({B_R})}}} \right|} ^p}dx + c{\int_{S({B_R})} {\left| {SXw} \right|} ^p}dx \\
  \le c{\left( {\frac{\rho }{R}} \right)^{Q + p}}{\int_{S({B_R})} {\left| {SXu - {{\left( {SXu} \right)}_{S({B_R})}}} \right|} ^p}dx + c{\int_{S({B_{2R}})} {\left| {SXu - {{\left( {SXu} \right)}_{S({B_{2R}})}}} \right|} ^p}dx \\
  \le \left( {c{{\left( {\frac{\rho }{R}} \right)}^{Q + p}} + c} \right){\int_{S({B_{2R}})} {\left| {SXu - {{\left( {SXu} \right)}_{S({B_{2R}})}}} \right|} ^p}dx. \\
 \end{array}$$
Now we use Lemma 2.8 to obtain (4.19).


\section{Proofs of Theorems 1.2, 1.3 and 1.4}

\label{5}

In order to prove Theorems 1.2, 1.3 and 1.4, we divide (1.1) into
two new systems and let $u=v+w$, where $v$ satisfies the nondiagonal
homogeneous system in $S(B_R )$
\begin{equation}
\label{eq37} \left\{ \begin{array}{l}
 - X_\alpha ^ *  \left( {a_{ij}^{\alpha \beta } X_\beta  v^j } \right) = 0, \\
v - u \in W_{X,0}^{1,2} \left( {S(B_R ),\mathbb{R}^N } \right) \\
 \end{array} \right.
\end{equation}
and $w $ solves the nondiagonal nonhomogeneous system in $S(B_R )$
\begin{equation}
\label{eq38} \left\{ \begin{array}{l}
 - X_\alpha ^ *  \left( {a_{ij}^{\alpha \beta } X_\beta  w^j } \right) = g_i(x,u,Xu)  - X_\alpha ^ *  f_i^\alpha (x,u,Xu), \\
w \in W_{X,0}^{1,2} \left( {S(B_R ),\mathbb{R}^N } \right). \\
 \end{array} \right.
\end{equation}
\noindent

\textbf{Proof of Theorem 1.2:} By (3.9) replacing $u$ by $w$ and
Theorem 4.1 replacing $u$ by $v$, it follows
$$\int_{S(B_R )} {\left| {SXw} \right|^p dx}  \le cR^\lambda
\left( {\left\| {Sg} \right\|_{L_X^{p,\lambda } }^p  + \left\|
{S\tilde {\tilde g}} \right\|_{L_X^{p,\lambda } }^p  + \left\|
{S\tilde g} \right\|_{L_X^{pq_0 ,\lambda q_0 } }^p } \right),$$
$$\int_{S(B_\rho  )} {\left| {SXv} \right|} ^p dx \le c\left(
{\frac{\rho } {R}} \right)^{\frac{{2Q - p(Q - \mu _1 )}} {2}}
\int_{S(B_R )} {\left| {SXv} \right|} ^p dx.$$ Using $u=v+w$, it
shows
\[
 \begin{gathered}
\quad  \int_{S(B_\rho  )} {\left| {SXu} \right|^p dx}  \le c\int_{S(B_\rho  )} {\left| {SXv} \right|^p dx}  + c\int_{S(B_\rho  )} {\left| {SXw} \right|^p dx}  \hfill \\
   \le c\left( {\frac{\rho }
{R}} \right)^{\frac{{2Q - p(Q - \mu _1 )}}
{2}} \int_{S(B_R )} {\left| {SXu} \right|} ^p dx + c\int_{S(B_R )} {\left| {SXw} \right|^p dx}  \hfill \\
   \le c\left( {\frac{\rho }
{R}} \right)^{\frac{{2Q - p(Q - \mu _1 )}}
{2}} \int_{S(B_R )} {\left| {SXu} \right|} ^p dx
+ cR^\lambda  \left( {\left\| {Sg} \right\|_{L_X^{p,\lambda } }^p
+ \left\| {S\tilde {\tilde g}} \right\|_{L_X^{p,\lambda } }^p
+ \left\| {S\tilde g} \right\|_{L_X^{pq_0 ,\lambda q_0 } }^p } \right). \hfill \\
\end{gathered}
\]
Taking $ H(\rho ) = \int_{S(B_\rho  )} {\left| {SXu} \right|^p dx},
H(R) = \int_{S(B_R )} {\left| {SXu} \right|^p dx}, a = \frac{{2Q -
p(Q - \mu _1 )}} {2}, b = \lambda, B = c\left( {\left\| {Sg}
\right\|_{L_X^{p,\lambda } }^p  + \left\| {S\tilde {\tilde g}}
\right\|_{L_X^{p,\lambda } }^p  + \left\| {S\tilde g}
\right\|_{L_X^{pq_0 ,\lambda q_0 } }^p } \right).$ Then there exists
$\mu _1$, $\frac{{\left( {p - 2} \right)Q + 2\lambda }} {p} < \mu _1
< Q $, such that $a>b$. We have by Lemma 2.8 that
\[
\int_{S(B_\rho  )} {\left| {SXu} \right|^p dx}  \le c\left(
{\frac{\rho } {R}} \right)^\lambda  \int_{S(B_R )} {\left| {SXu}
\right|} ^p dx + c\rho ^\lambda  \left( {\left\| {Sg}
\right\|_{L_X^{p,\lambda } }^p  + \left\| {S\tilde {\tilde g}}
\right\|_{L_X^{p,\lambda } }^p  + \left\| {S\tilde g}
\right\|_{L_X^{pq_0 ,\lambda q_0 } }^p } \right)
\]
and $SXu \in L_X^{p,\lambda } \left( {S(B_\rho  ),\mathbb{R}^N }
\right).$ Hence the result is proved.

\textbf{Proof of Theorem 1.3:} By Theorem 1.2, we see
\[
\int_{S(B_\rho  )} {\left| {SXu} \right|^p dx}  \le c\rho ^\lambda .
\]
Since $Q - p < \lambda  < Q$, it follows by taking $\kappa  = 1 -
\frac{{Q - \lambda }} {p}$ and using Lemma 2.10 that the conclusion
is true.

\textbf{Remark 5.1} Of course, we can also obtain H\"{o}lder
regularity by the isomorphic relationship between the Campanato
space $\mathcal {L}_{X}^{p,\lambda' }$($ - p < \lambda ' < 0$) and
the H\"{o}lder space $C_X^{0,\alpha }$($\alpha  =  - \frac{{\lambda
'}} {p}$) given in [8, Theorem 2.2].

Before proving Theorem 1.4, we first establish the following lemma.

\textbf{Lemma 5.1} Let $ w \in W_{X,0}^{1,2} (\Omega ,\mathbb{R}^N
)$ be a weak solution to (5.2), the coefficients $a_{ij}^{\alpha
\beta }$ in (5.2) satisfy (H1), $g_i (x,u,Xu)$ and $f_i^\alpha
(x,u,Xu)$ satisfy (H3). Then for any $p \in \left[ {2,2 +
\frac{{2Q}} {{Q + 2}}\varepsilon _0 } \right)$, $B_{2R}  \subset
\subset \Omega$, we have
\begin{eqnarray}
&& \quad \int_{B_R } {\left| {Xw} \right|^p dx} \notag \\
&& \le c\int_{B_{2R} } {\left| {Xu - (Xu)_{B_{2R} } } \right|} ^p dx
+ cR^\lambda  \left( {\left\| g \right\|_{L_X^{p,\lambda } (\Omega
)}^p  + \left\| {\tilde {\tilde g}} \right\|_{L_X^{p,\lambda }
(\Omega )}^p  + \left\| {\tilde g} \right\|_{L_X^{pq_0 ,\lambda q_0
} (\Omega )}^p } \right). \label{5.3}
\end{eqnarray}

\textbf{Proof}: Let us note that $w$ is also a weak solution to the
system
\[
 - X_\alpha ^ *  \left( {a_{ij}^{\alpha \beta } X_\beta  w^j } \right) = g_i  - X_\alpha ^ *  \left( {f_i^\alpha   - \left( {f_i^\alpha  } \right)_{B_{2R} } } \right).
\]
Multiplying both sides of the system by $w$ and integrating on
$B_{2R}$,
\begin{eqnarray}
&& \quad - \int_{B_{2R} } {A^{\alpha \beta } (x)\delta _{ij}
X_\alpha w^i X_\beta  w^j } dx \notag \\
&& = \int_{B_{2R} } {B_{ij}^{\alpha \beta } X_\alpha  w^i X_\beta
w^j } dx + \int_{B_{2R} } {g_i w^i } dx - \int_{B_{2R} } {\left(
{f_i^\alpha   - \left( {f_i^\alpha  } \right)_{B_{2R} } }
\right)X_\alpha  w^i } dx. \label{5.4}
\end{eqnarray}
By (H3), (2.3), H\"{o}lder's inequality and Young's inequality, it
implies
\[\begin{array}{l}
 \quad \int_{{B_{2R}}} {\left| {{g_i}} \right|\left| {{w^i}} \right|dx}  \le \int_{S({B_{2R}})} {\left( {{g^i} + L{{\left| {Xu} \right|}^{{\gamma _0}}}} \right)\left| {{w^i}} \right|dx}  \\
  \le {\left( {\int_{{B_{2R}}} {{{\left| {\tilde g} \right|}^{\frac{{2Q}}{{Q + 2}}}}dx} } \right)^{\frac{{Q + 2}}{{2Q}}}}{\left( {\int_{{B_{2R}}} {{{\left| w \right|}^{\frac{{2Q}}{{Q - 2}}}}dx} } \right)^{\frac{{Q - 2}}{{2Q}}}} \\
 {\kern 1pt} {\kern 1pt} {\kern 1pt} {\kern 1pt} {\kern 1pt} {\kern 1pt} {\kern 1pt}  + c{\left( {\int_{{B_{2R}}} {{{\left| {Xu} \right|}^2}dx} } \right)^{\frac{{{\gamma _0}}}{2}}}{\left( {\int_{{B_{2R}}} {{{\left| w \right|}^{\frac{2}{{2 - {\gamma _0}}}}}dx} } \right)^{\frac{{2 - {\gamma _0}}}{2}}} \\
  \le {\left( {\int_{{B_{2R}}} {{{\left| {\tilde g} \right|}^{2{q_0}}}dx} } \right)^{\frac{1}{{2{q_0}}}}}{\left( {\int_{{B_{2R}}} {{{\left| {Xw} \right|}^2}dx} } \right)^{\frac{1}{2}}} \\
 {\kern 1pt} {\kern 1pt} {\kern 1pt} {\kern 1pt} {\kern 1pt} {\kern 1pt} {\kern 1pt}  + c{\left( {\int_{{B_{2R}}} {{{\left| {Xu} \right|}^2}dx} } \right)^{\frac{{{\gamma _0}}}{2}}}R{\left| {{B_{2R}}} \right|^{\frac{{2 - {\gamma _0}}}{2}}}{\left( {\frac{1}{{\left| {{B_{2R}}} \right|}}\int_{{B_{2R}}} {{{\left| {Xw} \right|}^2}dx} } \right)^{\frac{1}{2}}} \\
  \le {c_\varepsilon }{\left( {\int_{{B_{2R}}} {{{\left| {\tilde g} \right|}^{2{q_0}}}dx} } \right)^{\frac{1}{{{q_0}}}}} + \varepsilon \int_{{B_{2R}}} {{{\left| {Xw} \right|}^2}dx}  + {c_\varepsilon }{\left( {\int_{{B_{2R}}} {{{\left| {Xu} \right|}^2}dx} } \right)^{{\gamma _0}}} \\
 {\kern 1pt} {\kern 1pt} {\kern 1pt} {\kern 1pt} {\kern 1pt} {\kern 1pt} {\kern 1pt} {\kern 1pt}  + \varepsilon {R^{Q + 2 - Q{\gamma _0}}}\int_{{B_{2R}}} {{{\left| {Xw} \right|}^2}dx} . \\
 \end{array}\]
Inserting it into (5.4), and using (H1) and (3.11), we have
\[\begin{array}{l}
 \quad {\lambda _0}\int_{{B_{2R}}} {{{\left| {Xw} \right|}^2}} dx \\
  \le \delta {\lambda _0}\int_{{B_{2R}}} {{{\left| {Xw} \right|}^2}} dx + \int_{{B_{2R}}} {\left| {{g_i}} \right|\left| {{w^i}} \right|} dx + \int_{{B_{2R}}} {\left| {f_i^\alpha  - {{\left( {f_i^\alpha } \right)}_{{B_{2R}}}}} \right|\left| {{X_\alpha }{w^i}} \right|} dx \\
  \le \delta {\lambda _0}\int_{{B_{2R}}} {{{\left| {Xw} \right|}^2}} dx + {c_\varepsilon }{\left( {\int_{{B_{2R}}} {{{\left| {\tilde g} \right|}^{2{q_0}}}dx} } \right)^{\frac{1}{{{q_0}}}}} + \varepsilon \int_{{B_{2R}}} {{{\left| {Xw} \right|}^2}dx}  \\
\quad + {c_\varepsilon }{\left( {\int_{{B_{2R}}} {{{\left| {Xu} \right|}^2}dx} } \right)^{{\gamma _0}}} + \varepsilon {R^{Q + 2 - Q{\gamma _0}}}\int_{{B_{2R}}} {{{\left| {Xw} \right|}^2}dx}  \\
\quad + {c_\varepsilon }\int_{{B_{2R}}} {{{\left| {f - {{\left( f \right)}_{{B_{2R}}}}} \right|}^2}} dx + \varepsilon \int_{{B_{2R}}} {{{\left| {Xw} \right|}^2}dx}  \\
  \le {c_\varepsilon }{\left( {\int_{{B_{2R}}} {{{\left| {\tilde g} \right|}^{2{q_0}}}dx} } \right)^{\frac{1}{{{q_0}}}}} + \left( {\delta {\lambda _0} + 2\varepsilon  + \varepsilon {R^{Q + 2 - Q{\gamma _0}}}} \right)\int_{{B_{2R}}} {{{\left| {Xw} \right|}^2}dx}  \\
\quad + {c_\varepsilon }{\left( {\int_{{B_{2R}}} {{{\left| {Xu} \right|}^2}dx} } \right)^{{\gamma _0}}} + {c_\varepsilon }\int_{{B_{2R}}} {{{\left| {\tilde {\tilde g} - {{\left( {\tilde {\tilde g}} \right)}_{{B_{2R}}}}} \right|}^2}dx}  + {c_\varepsilon }\int_{{B_{2R}}} {{{\left| {Xu - {{\left( {Xu} \right)}_{{B_{2R}}}}} \right|}^2}dx}  \\
  \le {c_\varepsilon }{\left( {\int_{{B_{2R}}} {{{\left| {\tilde g} \right|}^{2{q_0}}}dx} } \right)^{\frac{1}{{{q_0}}}}} + {\theta _3}\int_{{B_{2R}}} {{{\left| {Xw} \right|}^2}dx}  + {c_\varepsilon }{\left( {\int_{{B_{2R}}} {{{\left| {Xu} \right|}^2}dx} } \right)^{{\gamma _0}}} \\
\quad  + {c_\varepsilon }\int_{{B_{2R}}} {{{\left| {\tilde {\tilde g}} \right|}^2}dx}  + {c_\varepsilon }\int_{{B_{2R}}} {{{\left| {Xu - {{\left( {Xu} \right)}_{{B_{2R}}}}} \right|}^2}dx} , \\
 \end{array}\]
where $\theta _3  = \delta \lambda _0  + 2\varepsilon  + \varepsilon
R^{Q + 2 - Q\gamma _0 }$. By choosing $\varepsilon$ small enough
such that $ \lambda _0  - \theta _3  > 0$ and applying (3.8), we
obtain
\begin{eqnarray}
&& \quad \int_{B_{2R} } {\left| {Xw} \right|^2 } dx \notag \\
 && \le c_\varepsilon  \left( {\int_{B_{2R} } {\left| {\tilde g} \right|^{2q_0 } dx} } \right)^{\frac{1}
{{q_0 }}}  + c_\varepsilon  \int_{B_{2R} } {\left| {Xu} \right|^2
dx} + c_\varepsilon  \int_{B_{2R} } {\left| {\tilde {\tilde g}}
\right|^2 dx} \notag \\
&& \quad + c_\varepsilon  \int_{B_{2R} } {\left| {Xu - \left( {Xu} \right)_{B_{2R} } } \right|^2 dx}  \notag \\
&& \le c\left( {\int_{B_{2R} } {\left| {\tilde g} \right|^{2q_0 }
dx} } \right)^{\frac{1} {{q_0 }}}  + c\int_{B_{2R} } {\left| g
\right|^2 dx} + c\int_{B_{2R} } {\left| {\tilde {\tilde g}}
\right|^2 dx}
+ c\int_{B_{2R} } {\left| {Xu - \left( {Xu} \right)_{B_{2R} } } \right|^2 dx}  \notag \\
&& \le cR^{\frac{{2\lambda }} {p}} \left| {B_{2R} }
\right|^{\frac{{p - 2}} {{pq_0 }}} \left\| {\tilde g}
\right\|_{L_X^{pq_0 ,\lambda q_0 } }^2  + cR^{\frac{{2\lambda }}
{p}} \left| {B_{2R} } \right|^{\frac{{p - 2}} {p}} \left( {\left\| g
\right\|_{L_X^{p,\lambda } }^2
+ \left\| {\tilde {\tilde g} } \right\|_{L_X^{p,\lambda } }^2 } \right) \notag \\
&& \quad + c\left| {B_{2R} } \right|^{\frac{{p - 2}} {p}} \left(
{\int_{B_{2R} } {\left| {Xu - \left( {Xu} \right)_{B_{2R} } }
\right|^p dx} } \right)^{\frac{2} {p}}. \label{5.5}
\end{eqnarray}
From Lemma 3.2, it yields
\[
\begin{gathered}
  \int_{B_R } {\left| {Xw} \right|^p dx}  \le c\left| {B_R } \right|^{\frac{{2 - p}}
{2}} \left( {\int_{B_{2R} } {\left| {Xw} \right|^2 dx} }
\right)^{\frac{p}
{2}}  \hfill \\
 \quad \quad  \quad \quad \quad \quad + c\int_{B_{2R} } {\left( {\left| g \right|^p
 + \left| {\tilde {\tilde g}} \right|^p } \right)dx}  + c\left| {B_R } \right|^{\frac{{q_0  - 1}}
{{q_0 }}} R^p \left( {\int_{B_{2R} } {\left| {\tilde g}
\right|^{pq_0 } dx} } \right)^{\frac{1}
{{q_0 }}} . \hfill \\
\end{gathered}
\]
Putting (5.5) into it and noting (2.1), we have
\[
\begin{gathered}
  \quad \int_{B_R } {\left| {Xw} \right|^p dx}  \hfill \\
   \le cR^{p + \lambda  - 2} \left\| {\tilde g} \right\|_{L_X^{pq_0 ,\lambda q_0 } }^p
   + cR^\lambda  \left( {\left\| g \right\|_{L_X^{p,\lambda } }^p  + \left\| {\tilde {\tilde g}}
   \right\|_{L_X^{p,\lambda } }^p } \right) + c\int_{B_{2R} } {\left| {Xu - \left( {Xu} \right)_{B_{2R} } } \right|^p dx}  \hfill \\
   \le cR^\lambda  \left( {\left\| {\tilde g} \right\|_{L_X^{pq_0 ,\lambda q_0 } }^p
   + \left\| g \right\|_{L_X^{p,\lambda } }^p  + \left\| {\tilde{ \tilde g}} \right\|_{L_X^{p,\lambda } }^p } \right) + c\int_{B_{2R} } {\left| {Xu - \left( {Xu} \right)_{B_{2R} } } \right|^p dx} . \hfill \\
\end{gathered}
\]
It completes the proof.

\textbf{Proof of Theorem 1.4: } By Lemma 5.1, it follows
\[\begin{array}{l}
 \quad \int_{S({B_R})} {{{\left| {SXw} \right|}^p}dx}  \\
  \le c{R^\lambda }\left( {\left\| {Sg} \right\|_{L_X^{p,\lambda }}^p + \left\| {S\tilde {\tilde g}} \right\|_{L_X^{p,\lambda }}^p + \left\| {S\tilde g} \right\|_{L_X^{p{q_0},\lambda {q_0}}}^p} \right) + c\int_{S({B_{2R}})} {{{\left| {SXu - {{\left( {SXu} \right)}_{S({B_{2R}})}}} \right|}^p}dx} . \\
 \end{array}\]
Using $u=v+w$, Theorem 4.2 and the above inequality, we have
\[\begin{array}{l}
 \quad \int_{S({B_\rho })} {{{\left| {SXu - {{\left( {SXu} \right)}_{S({B_\rho })}}} \right|}^p}dx}  \\
  \le c\int_{S({B_\rho })} {{{\left| {SXv - {{\left( {SXv} \right)}_{S({B_\rho })}}} \right|}^p}dx}
  + c\int_{S({B_\rho })} {{{\left| {SXw - {{\left( {SXw} \right)}_{S({B_\rho })}}} \right|}^p}dx}  \\
  \le c{\left( {\frac{\rho }{R}} \right)^{{\mu _2}}}\int_{S({B_R})} {{{\left| {SXv - {{\left( {SXv} \right)}_{S({B_R})}}} \right|}^p}dx}
  + c\int_{S({B_\rho })} {{{\left| {SXw} \right|}^p}dx}  \\
  \le c{\left( {\frac{\rho }{R}} \right)^{{\mu _2}}}\int_{S({B_R})} {{{\left| {SXu - {{\left( {SXu} \right)}_{S({B_R})}}} \right|}^p}dx}
  + c\int_{S({B_R})} {{{\left| {SXw} \right|}^p}dx}  \\
  \le c{\left( {\frac{\rho }{R}} \right)^{{\mu _2}}}\int_{S({B_{2R}})} {{{\left| {SXu - {{\left( {SXu} \right)}_{S({B_{2R}})}}} \right|}^p}dx}  \\
 {\kern 1pt} {\kern 1pt} {\kern 1pt} {\kern 1pt} {\kern 1pt} {\kern 1pt} {\kern 1pt} {\kern 1pt}
 + c{R^\lambda }\left( {\left\| {S\tilde g} \right\|_{L_X^{p{q_0},\lambda {q_0}}}^p + \left\| {Sg} \right\|_{L_X^{p,\lambda }}^p
 + \left\| {S\tilde {\tilde g}} \right\|_{L_X^{p,\lambda }}^p} \right)
 + c\int_{S({B_{2R}})} {{{\left| {SXu - {{\left( {SXu} \right)}_{S({B_{2R}})}}} \right|}^p}dx}  \\
  \le \left[ {c{{\left( {\frac{\rho }{R}} \right)}^{{\mu _2}}} + c} \right]\int_{S({B_{2R}})} {{{\left| {SXu - {{\left( {SXu} \right)}_{S({B_{2R}})}}} \right|}^p}dx}  \\
 {\kern 1pt} {\kern 1pt} {\kern 1pt} {\kern 1pt} {\kern 1pt} {\kern 1pt} {\kern 1pt}
 + c{R^\lambda }\left( {\left\| {S\tilde g} \right\|_{L_X^{p{q_0},\lambda {q_0}}}^p + \left\| {Sg} \right\|_{L_X^{p,\lambda }}^p
 + \left\| {S\tilde {\tilde g}} \right\|_{L_X^{p,\lambda }}^p} \right). \\
 \end{array}\]
Set $H(\rho ) = \int_{S(B_\rho  )} {\left| {SXu - \left( {SXu}
\right)_{S(B_\rho  )} } \right|^p dx}, H(R) = \int_{S(B_{2R} )}
{\left| {SXu - \left( {SXu} \right)_{S(B_{2R} )} } \right|^p dx}, a
= \mu _2, b = \lambda, B = c\left( {\left\| {Sg}
\right\|_{L_X^{p,\lambda } }^p  + \left\| {S\tilde {\tilde g}}
\right\|_{L_X^{p,\lambda } }^p  + \left\| {S\tilde g}
\right\|_{L_X^{pq_0 ,\lambda q_0 } }^p } \right)$. Using $0 < \mu _2
< Q + p, 0 < \lambda  < Q$, it derives that there exists $\mu _2$
such that $\mu _2 > \lambda $. We have by Lemma 2.8 that
\[
\begin{gathered}
 \quad \int_{S(B_\rho  )} {\left| {SXu - \left( {SXu} \right)_{S(B_\rho  )} } \right|^p dx}  \hfill \\
   \le c\left( {\frac{\rho }
{R}} \right)^\lambda  \int_{S(B_{2R} )} {\left| {SXu - \left( {SXu} \right)_{S(B_{2R} )} } \right|^p dx}
  + c\rho ^\lambda  \left( {\left\| {S\tilde g} \right\|_{L_X^{pq_0 ,\lambda q_0 } }^p
  + \left\| {Sg} \right\|_{L_X^{p,\lambda } }^p  + \left\| {S\tilde {\tilde g}} \right\|_{L_X^{p,\lambda } }^p } \right). \hfill \\
\end{gathered}
\]
Hence \[ SXu \in \mathcal {L}_X^{p,\lambda } (S(B_\rho
),\mathbb{R}^N )
\]
and the proof is finished.


{\bf Acknowledgements } This work is supported by the National
Natural Science Foundation of China(Grant Nos. 11271299, 11001221);
Natural Science Foundation Research Project of Shaanxi Province
(Grant No. JC201124).


\end{document}